\documentclass[12pt,envcountsame]{article}
\usepackage{amssymb}
\usepackage{amsmath,amsthm}
\usepackage[pdfborder={0 0 0}]{hyperref}
\newtheorem{theorem}{Theorem}[section]

\newtheorem{lemma}[theorem]{Lemma}
\newtheorem{definition}[theorem]{Definition}
\newtheorem{corollary}[theorem]{Corollary}
\newtheorem{proposition}[theorem]{Proposition}
\newtheorem{problem}[theorem]{Problem}
 
\def\RR{{\mathbb R}}
\def\OL{\overline}

\def\({\left(}
\def\){\right)}
\def\[{\left[}
\def\]{\right]}

\def \W{\widetilde}

\begin {document}

\title {Pointwise optimal multivariate spline method \\ for recovery of twice differentiable functions on a simplex}

\author {Sergiy V. Borodachov \\
\\
{\it Department of Mathematics, Towson University,}\\
{\it Towson, MD, 21252, USA}\\
\\
e-mail: sborodachov@towson.edu}

\maketitle

\begin {abstract}
We obtain the spline recovery method on a $d$-dimensional simplex $T$ that uses as information values and gradients of a function $f$ at the vertices of $T$ and is optimal for recovery of $f({\bf w})$ at every point ${\bf w}$ of an admissible domain $P\supset T$ on the class $W^2(P)$ of twice differentiable functions on $P$ with uniformly bounded second order derivatives in any direction. If, in particular, every face of $T$ (of any dimension) contains its circumcenter, we can take $P=T$. 

We also find the error function of the pointwise optimal method which turns out to be a function in $W^2(P)$ with zero information. The error function is a piecewise quadratic $C^1$-function over a certain polyhedral partition and can be considered as a multivariate analogue of the classical Euler spline $\phi_2$. The pointwise optimal method is a continuous spline of degree two (with some pieces of degree one) over the same partition.

\end {abstract}

\section {Rigorous setting of the problem}

Let $T\subset \RR^d$ be a non-degenerate simplex with vertices ${\bf v}_0,{\bf v}_1,\ldots,{\bf v}_d$; i.e., a simplex that is not contained in any $(d-1)$-dimensional hyperplane.

Let $\mathcal F\subset C^1(T)$ be a class of functions (that will be specified later). Here $C^1(T)$ denotes the space of continuously differentiable functions on $T$. Assume that for every $f\in \mathcal F$, the following information vector
$$
I(f):=\(f({\bf v}_0),\ldots,f({\bf v}_d),\nabla f({\bf v}_0),\ldots,\nabla f ({\bf v}_d)\)\in \RR^{(d+1)^2}
$$
is known. We pursue the following two goals:

(i) based on the data given by the vector $I(f)$, recover the value of a function $f\in \mathcal F$ at a given point ${\bf w}\in T$ (in an optimal way);

(ii) based on the data given by the vector $I(f)$, produce (in an optimal way) a continuous function on $T$ that approximates the function $f\in \mathcal F$ on the whole domain $T$ (global recovery).

The worst-case error setting of the optimal recovery problem considered in this paper goes back to the works by A.N.~Kolmogorov in the 1940s. First results on this problem were obtained by A.~Sard \cite {Sar49}, S.M.~Nikol'skii \cite {Nik50}, and J.~Kiefer \cite {Kif57}. Later, a large number of optimal recovery results were obtained by many authors. Detailed reviews of these results can be found in the books \cite {TraWozOTOA}, \cite {KorSTP}, \cite {TraWozWasIBC}, \cite {KorTCTP}, \cite {MotLygDor}, \cite {Zhe03}, and \cite {NovWozTMP}.

Due to technicalities in our arguments we have to assume that $\mathcal F\subset C^1(P)$, where
$P$ is some closed convex set in $\RR^d$ that contains $T$ (we will specify additional assumptions on $P$ in Definition \ref {D3.1}). We will call the set $P$ an admissible domain.
If $T$ is completely well-centered (see Definition \ref {cwc}), we will be able to assume that $P=T$. 

Since our arguments stay the same, we can assume that ${\bf w}\in P$ in goal (i) and recover $f$ on the set $P$ rather than $T$ in goal (ii).

To achieve goal (i), we consider the family of all non-adaptive recovery algorithms that use the information vector $I(f)$ and look for optimal ones on the class $\mathcal F$.
From the most general point of view, any non-adaptive algorithm for the recovery of the value $f({\bf w})$ that uses the data vector $I(f)$ as an input, produces some number as an output; i.e., it is defined by some mapping $\Psi:\RR^{(d+1)^2}\to\RR$ via the formula
\begin {equation}\label {PSI}
s(f;\Psi):=\Psi\(I(f)\).
\end {equation}
Conversely, every mapping $\Psi:\RR^{(d+1)^2}\to\RR$ defines a non-adaptive algorithm of form \eqref {PSI} for recovering the value $f({\bf w})$. The error of any algorithm \eqref {PSI} over the class $\mathcal F$ is defined by
$$
e_{\bf w}(\mathcal F;\Psi):=\sup\limits_{f\in \mathcal F}{\left|f({\bf w})-\Psi(I(f))\right|}.
$$
\begin {problem}\label {P1}
{\rm
Given a point ${\bf w}\in P$ and a class $\mathcal F\subset C^1(P)$, find the value
\begin {equation}\label {e1}
e_{\bf w}(\mathcal F):=\inf\limits_{\Psi:\RR^{(d+1)^2}\to \RR}{e_{\bf w}(\mathcal F;\Psi)}
\end {equation}
and optimal algorithms for the recovery of the value $f({\bf w})$ on the class $\mathcal F$; i.e., mappings $\Psi^\ast$ that attain the infimum on the right-hand side of \eqref {e1} (if they exist).
}
\end {problem}

To state rigorously the problem of optimal global recovery mentioned in the goal (ii) we will need the following notation. Let $C(P)=C^0(P)$ be the space of functions continuous on the set $P$ and let $\|\cdot\|$ be a monotone seminorm on $C(P)$; i.e., a seminorm such that $\|f_1\|\leq \|f_2\|$ whenever the inequality  $\left|f_1({\bf x})\right|\leq \left|f_2({\bf x})\right|$ holds at every point ${\bf x}\in P$. This is a natural assumption on the seminorm: if the absolute value of the difference between the function and its approximant increases at every point, the value of the approximation error also increases. Important examples of a monotone seminorm are the uniform norm on $P$, the standard $L_p$-norm on $P$, $1\leq p<\infty$, as well as the weighted $L_p$-norm on $P$, which are often used in mathematical literature to measure the approximation error. Observe also that a monotone seminorm on $C(T)$ is a monotone seminorm on $C(P)$.

In the case of global recovery we will include into our considerations recovery algorithms that produce continuous approximants even though the functions being recovered are $C^1$-smooth. This allows to consider a broader set of algorithms and include the pointwise optimal one that turns out to be only continuous. Every mapping $\Phi:\RR^{(d+1)^2}\to C(P)$ defines a non-adaptive algorithm for global recovery of functions from the class $\mathcal F$ of the form
\begin {equation}\label {alg2}
S(f;\Phi)({\bf x}):=\Phi(I(f))({\bf x}),\ \ \ {\bf x}\in P.
\end {equation}
The error of method \eqref {alg2} over the class $\mathcal F$ in a given seminorm $\|\cdot \|$ is defined by
$$
R(\mathcal F;\Phi,\|\cdot\|):=\sup\limits_{f\in \mathcal F}{\|f-\Phi(I(f))\|}.
$$
\begin {problem}\label {P2}
{\rm
Given a class $\mathcal F\subset C^1(P)$ and a monotone seminorm $\|\cdot\|$ on $C(P)$, find the value
\begin {equation}\label {e2}
R(\mathcal F;\|\cdot\|):=\inf\limits_{\Phi:\RR^{(d+1)^2}\to C(P)}{R(\mathcal F;\Phi,\|\cdot \|)}
\end {equation}
and optimal algorithms for global recovery of functions from the class $\mathcal F$ in the seminorm $\|\cdot\|$; i.e., mappings $\Phi^\ast$ that attain the infimum on the right-hand side of \eqref {e2} (if they exist).
}
\end {problem}

\begin {definition}
{\rm
A global recovery algorithm $S(f;\Phi)$ of form \eqref {alg2} is called {\it pointwise optimal on a class $\mathcal F\subset C^1(P)$} if its restriction $s(f):=S(f;\Phi)({\bf w})$ to any point ${\bf w}\in P$ is optimal on the class $\mathcal F$ in the sense of Problem \ref {P1} for the recovery of the value $f({\bf w})$.
}
\end {definition}
If $S(f;\Phi)$ is an algorithm of form \eqref {alg2} optimal on the class $\mathcal F$ in any monotone seminorm $\|\cdot\|$, letting $\|f\|=\left|f({\bf w})\right|$ where ${\bf w}$ is arbitrary fixed point in $P$, and restricting $S(f;\Phi)$ just to the point ${\bf w}$, we obtain an optimal algorithm (in the sense of Problem \ref {P1}) for recovery of $f({\bf w})$ on the class $\mathcal F$; i.e., $S(f;\Phi)$ is pointwise optimal.

We conclude the introduction by describing the class $\mathcal F$ on which Problems \ref {P1} and \ref {P2} will be considered. Many works on optimal recovery study the one-dimensional class $W^r_\infty[a,b]$ of functions $f:[a,b]\to \RR$ such that $f^{(r-1)}$ is absolutely continuous on $[a,b]$ and $\left|f^{(r)}(t)\right|\leq 1$ for a.e. $t\in [a,b]$ (we agree here that $f^{(0)}=f$ in the case $r=1$). A function $f:[a,b]\to\RR$ belongs to the class $W^r_\infty[a,b]$ if and only if $f^{(r-1)}$ is defined everywhere on $[a,b]$ and $\left|f^{(r-1)}(t_1)-f^{(r-1)}(t_2)\right|\leq \left|t_1-t_2\right|$ for all $t_1,t_2\in [a,b]$.

\begin {definition}
{\rm
Let $W^r(P)$, $r=1,2$, be the class of functions $f\in C^{r-1}(P)$ such that for any two points ${\bf x},{\bf y}\in P$, ${\bf x}\neq {\bf y}$, the function $g(t):=f\({\bf x}+\frac {t}{\rho}\({\bf y}-{\bf x}\)\)$, $t\in [0,\rho]$, where $\rho=\left|{\bf y}-{\bf x}\right|$, belongs to the class $W^r_\infty [0,\rho]$. 
}
\end {definition}Notice that $W^1(P)$ is the Lipschitz class of functions defined on $P$ and the pointwise optimal method of recovery from the data given by $f({\bf v}_0),\ldots,f({\bf v}_d)$ is known for this class, see \cite {Bab76a} (one can also see \cite {BabBojBor}). 
Our main interest is finding a pointwise optimal recovery method for functions from the class $W^2(P)$ based on the data vector $I(f)$ and finding the optimal error function $\W g({\bf w}):=e_{\bf w}(W^2(P))$, ${\bf w}\in P$. 

The problem of finding the error function $\W g({\bf w})$ was suggested to the author by V.F. Babenko. Since the class $W^2(P)$ is convex and centrally symmetric, the motivation for his question is the following. 
If $\W g$ is in the class $W^2(P)$ with $I(\W g)={\bf 0}$, and for every given $f\in W^2(P)$,
the function defined on $P$ by combining the output values of optimal algorithms for recovery of $f({\bf w})$ at every ${\bf w}\in P$ is continuous on $P$, {\it then} the resulting global recovery algorithm will be optimal on the class $ W^2(P)$ in any monotone seminorm $\|\cdot \|$ and the optimal recovery error will be given by $\|\W g\|$.

The paper is structured as follows. Section \ref {S2''} contains the review of known results on pointwise optimal recovery algorithms. Section \ref {S3add} discusses admissible domains for the technique we use. In particular, it states the admissibility of a completely well-centered simplex (see Proposition \ref {T_c}). Section \ref {S4results} states the pointwise optimal recovery method (Theorem~\ref {Th1} and its immediate consequence, Corollary \ref {wcT}). Section \ref {S4results} also describes certain properties of the error function $\mathcal E$ (Theorem \ref {Th2} and an immediate consequence, Corollary \ref {q1'}). Section \ref {S5p} discusses further properties of the error function $\mathcal E$ (Theorems \ref {E} and \ref {nd}) and of the underlying partition (Proposition \ref {partition} and Theorem \ref {well_c}). Section \ref {aux} contains proofs of several auxiliary statements. Section \ref {T_c1} proves Proposition \ref {T_c} and Theorem \ref {well_c}. Section \ref {S8_4.3} proves Theorem \ref {Th2}. Section \ref {S9_4.1} proves Theorem \ref {Th1} while Section \ref {S10} contains the proofs of the rest of the statements from Section~\ref {S5p}. Section \ref {A11} is the Appendix of this paper.

\section {Review of known results}\label {S2''}

Recovery problems \ref {P1} and \ref {P2} on classes of differentiable functions have been extensively studied in the one-dimensional case. Detailed reviews of these results are given in the books \cite {TraWozOTOA}, \cite {KorSTP}, \cite {TraWozWasIBC}, \cite {KorTCTP}, \cite {MotLygDor}, \cite {Zhe03}, and \cite {NovWozTMP}. We will mention in detail only those one-dimensional results that prove pointwise optimality. A pointwise optimal algorithm for recovery of $f(t)$ on the class $W^r_\infty [a,b]$ that uses function values and derivatives up to order $l_1$ at $a$ and up to order $l_2$ at $b$ was obtained by B.D. Bojanov in \cite {Boj1975} in the case $l_1=l_2=r-1$. We note that the pointwise optimal algorithm here is a piecewise polynomial function (spline) of degree $r-1$ that has $r-2$ continuous derivatives. In the case $0\leq l_1,l_2\leq r-1$, $l_1+l_2\geq r$, the error function was described in the paper \cite {MicRivWin1976} and the pointwise optimal spline was described in the book \cite {KorTCTP} (see Section 5.3 in this book and references on p. 250). 

In the multivariate case, for the class of Lipschitz functions defined on a compact set $Q\subset \RR^d$, the pointwise optimal recovery algorithm that uses function values at $n$ points ${\bf x}_1,\ldots,{\bf x}_n$ in $Q$ was essentially described by Babenko, see \cite {Bab76a} (see also, for example, \cite {BabBojBor}). The output of the optimal algorithm is a piecewise constant function, which is equal to $f({\bf x}_i)$ on the intersection of $Q$ with the Voronoi cell of the node ${\bf x}_i$, $i=1,\ldots,n$ (there is a certain freedom in defining the values of the output function of the optimal algorithm on the intersection of $Q$ with the boundaries of the Voronoi cells). The pointwise optimal recovery algorithm on the class $W^2(T)$ that uses {\it only function values} at the vertices ${\bf v}_0,\ldots,{\bf v}_d$ of the simplex $T$ was found in \cite {Klz1996}. The linear $d$-variate polynomial that interpolates the values $f({\bf v}_i)$, $i=0,1,\ldots,d$, turns out to be pointwise optimal.

The problem of global recovery (when $\left\|\ \!\cdot\ \!\right\|$ is the uniform norm) of functions from the Lipschitz class $W^1(Y)$, where $Y$ is a convex $d$-dimensional polytope,  based on known values on a fixed finite subset $X$ of $Y$ that includes the vertices of $Y$, was solved in \cite {BabLyg1975}, \cite {Bab1978}, and \cite {Bab1980}. The optimality of the linear interpolating spline over a Delaunay triangulation of $X$ in $Y$ has been established in these papers. The author of the current paper together with T.S. Sorokina (see \cite {BorSor2010}) constructed a quadratic polynomial recovery method on $T$ (that uses vector $I(f)$ as information), which is optimal in the sense of Problem \ref {P2} on the class $W^2(T)$ when $\|\cdot\|$ is the uniform norm. This method partially interpolates the data. When $d=2$ and $T$ is an equilateral triangle we also constructed two fully interpolating polynomial methods that use $I(f)$ as information and are optimal on the class $W^2(T)$ in the sense of Problem \ref {P2} (when $\|\cdot \|$ is the uniform norm), see \cite {BorSor2012}.

The current paper finds a pointwise optimal algorithm on the class $W^2(P)$ ($d>1$), where $P$ is an admissible domain for the simplex $T$ (see Definition \ref {D3.1}), that uses {\it both function values and gradients} at each vertex of $T$. When $T$ is completely well-centered, we have $P=T$. This method is also optimal in the sense of Problem \ref {P2} on the class $W^2(P)$ for global recovery in any monotone seminorm $\|\cdot \|$. This optimal method is a linear-quadratic spline over a certain polyhedral partition of $\RR^d$ (and hence of $P$) that interpolates both function values and gradients at the vertices of $T$. It is continuous but is not smooth in general. The optimal method represents a multivariate generalization of the pointwise optimal method found by Bojanov \cite {Boj1975} on the univariate class $W^2_\infty [a,b]$. In the univariate case, we have $T=P=[a,b]$ and if we let ${ c}=(a+b)/2$, $u_0=(a+{c})/2$ and $u_1=(b+{c})/2$, then the one-dimensional optimal method is given by
$$
\W s(f;t)=\begin {cases} f(a)+f'(a)(t-a), & t\in [a,u_0],\\
b_0(f(a)+f'(a)(u_0-a))+b_1(f(b)+f'(b)(u_1-b)), & t\in (u_0,u_1),\\
f(b)+f'(b)(t-b), & t\in [u_1,b].\\ \end {cases}
$$
where $b_0=\frac {u_1-t}{u_1-u_0}$ and $b_1=\frac {t-u_0}{u_1-u_0}$ are the barycentric coordinates of the point $t$ in $[u_0,u_1]$. This method is a piecewise linear continuous function that interpolates the data vector $I(f)=(f(a),f(b),f'(a),f'(b))$. For every fixed $t\in [a,b]$, the error over the class $W^2_\infty[a,b]$ of the method $\W s(f;t)$ for recovery of the value $f(t)$ is given by the classical Euler spline
$$
\phi_2(t)=\begin {cases}\frac {(t-a)^2}{2},& t\in [a,u_0],\\
\frac {R^2}{4}-\frac {(t-{c})^2}{2}, & t\in (u_0,u_1),\\
\frac {(t-b)^2}{2},& t\in [u_1,b],\\
\end {cases}
$$
where $R=\left|{c}-a\right|=\left|{c}-b\right|$. Notice that the function $\phi_2$ is contained in the class $W^2_\infty[a,b]$ and vanishes at points $a$ and $b$ together with its first derivatives, that is, $I(\phi_2)=(0,0,0,0)$. Furthermore, the derivative $\phi_2'$ is a piecewise linear function whose pieces have slope $1$ or $-1$; i.e. each piece of $\phi_2$ is an isometry.

\section {Admissible domains}\label {S3add}

When obtaining a sharp error estimate for the optimal method on the simplex $T$ we will use function values at a certain set of key points (different from the vertices of $T$) and thus, need to make sure that these points lie in the domain of the function being recovered. This is why we make additional assumptions about the simplex $T$ or extend the domain of the functions being recovered to a larger set $P$ that contains all the key points. Certain notation used throughout the rest of the paper is also introduced in this section.

Denote by ${\bf c}$ the circumcenter of $T$ and let
$$
U:=\frac {1}{2}\({\bf c}+T\). 
$$
Then $U$ is the simplex in $\RR^d$ with vertices ${\bf u}_i=\frac 12 \({\bf c}+{\bf v}_i\)$, $i=0,1,\ldots,d$.
Let also $b_i=b_i({\bf x})$, $i=0,1,\ldots,d$, be the barycentric coordinates of a point ${\bf x}\in U$ relative to the simplex $U$, that is, the unique vector of numbers $(b_0({\bf x}),\ldots,b_d({\bf x}))$ such that
$$
{\bf x}=\sum_{i=0}^{d}b_i({\bf x}){\bf u}_i,
$$
$\sum_{i=0}^{d}b_i({\bf x})=1$, and $b_i({\bf x})\geq 0$, $i=0,1,\ldots,d$. Observe that $b_i({\bf x})$ is a linear function of ${\bf x}$, $b_i({\bf u}_i)=1$, and $b_i({\bf u}_j)=0$, $j\neq i$.
Let $\Sigma:=\{0,1,\ldots,d\}$ and for every non-empty subset $I\subset \Sigma$, define the set
$$
G_I:=\{{\bf x}\in U : b_i({\bf x})>0,\    i\in I, \ \ {\rm and}\ \ b_i({\bf x})=0,\ i\notin I \}.
$$
Notice that the set $G_\Sigma$ is the interior of $U$, the sets $G_{\{i\}}$, $i=0,1,\ldots,d$, represent the vertices of $U$ and all other sets $G_I$ are the relative interiors of faces of $U$ of dimensions from $1$ to $d-1$. Furthermore, the sets $G_I$, $I\subset \Sigma$, $I\neq \emptyset$, are pairwise disjoint.

Since $U$ is a compact convex set and the Euclidean norm is strictly convex, for every ${\bf x}\in \RR^d$, there is a {\it unique} point ${\bf y}$ in $U$ that is closest to ${\bf x}$ among all points in $U$ (see Lemma \ref {distconv}). Let $\varphi:\RR^d\to U$ be the mapping defined by $\varphi({\bf x})={\bf y}$. This mapping will play a crucial role in our further arguments. 

\begin {definition}\label {D3.1}
{\rm
Given a non-degenerate simplex $T$ in $\RR^d$, we call a set $P\subset \RR^d$ an {\it admissible domain for $T$} if $P$ is a convex set that contains $T$, and for every ${\bf x}\in P$, we have ${\bf x}-\varphi({\bf x})+{\bf u}_i\in P$, $i\in I$, where $I\subset \Sigma$ is the subset such that $\varphi({\bf x})\in G_I$.
}
\end {definition}
Notice that $\RR^d$ is an admissible domain for any simplex $T$. In the case $d=1$, the interval $T=[a,b]$ is admissible for itself. In the case $d=2$, a triangle is an admissible domain for itself if and only if it is non-obtuse; i.e., if and only if it contains its circumcenter. However, when $d= 3$, a simplex containing its circumcenter need not be an admissible domain for itself. Consider, for example, the simplex $T$ with vertices ${\bf v}_0(0,-1,0)$, ${\bf v}_1(\sqrt{1-a^2},a,0)$, ${\bf v}_2(-\sqrt{1-a^2},a,0)$, and ${\bf v}_3(0,b,\sqrt{1-b^2})$, where $0<a<b/2<1/2$. The circumcenter ${\bf c}(0,0,0)$ of $T$ belongs to $T$. The orthogonal projection ${\bf p}_3$ of the point ${\bf u}_3=({\bf c}+{\bf v}_3)/2$ onto the plane $\Pi$ passing through the points ${\bf v}_1$, ${\bf v}_2$, ${\bf v}_3$ belongs to the interior of the facet ${\bf v}_1{\bf v}_2{\bf v}_3$ of $T$. However, for $i=1,2$, the orthogonal projection ${\bf p}_i$ of the point ${\bf u}_i=({\bf c}+{\bf v}_i)/2$ onto the plane $\Pi$ does not belong to the facet ${\bf v}_1{\bf v}_2{\bf v}_3$. Then for any point ${\bf x}\in {\bf v}_1{\bf v}_2{\bf v}_3$ that also belongs to the interior of the triangle ${\bf p}_1{\bf p}_2{\bf p}_3$, we have $\varphi({\bf x})\in G_{\{1,2,3\}}$ and ${\bf x}-\varphi({\bf x})+{\bf u}_i={\bf p}_i\notin T$, $i=1,2$. Consequently, $T$ is not admissible for itself unless we make a stronger assumption on $T$.

Recall that a set of points $\{{\bf v}_i\}_{i=0}^k$ in $\RR^d$ is called {\it affinely independent}, if $\sum_{i=0}^{k}c_i{\bf v}_i\neq {\bf 0}$ for any non-zero vector of coefficients $(c_0,\ldots,c_k)$ such that $\sum_{i=0}^{k}c_i=0$.

\begin {definition}
{\rm
Given an affinely independent set of points $\{{\bf v}_i\}_{i\in J}$ in $\RR^d$,  denote by
$Q_J$ 
the affine subspace of $\RR^d$ spanned by the point set $\{{\bf v}_i\}_{i\in J}$; i.e.,
$
Q_J:=\left\{\sum_{i\in J}c_i{\bf v}_i : \sum_{i\in J}c_i=1\right\}.
$
Denote by $T_J$ the simplex with the set of vertices $\{{\bf v}_i\}_{i\in J}$. Let $H_J$ be the set of all points in $\RR^d$ that are equidistant from the vertices ${\bf v}_i$, $i\in J$. Denote also by ${\bf c}_J$ the circumcenter of $T_J$; i.e., the point in $Q_J$ equidistant from all the vertices ${\bf v}_i$, $i\in J$, and let $R_J$ be the circumradius of $T_J$. 
}
\end {definition}
\begin {definition}\label {cwc}{\rm
A simplex $T\subset \RR^d$ is called {\it completely well-centered} if $T$ and every face of $T$ (of any dimension) contain their circumcenters.
}
\end {definition}

We next describe a class of simplices in $\RR^d$ which are admissible domains. 

\begin {proposition}\label {T_c}
Any non-degenerate completely well-centered simplex $T$ in $\RR^d$ is an admissible domain for itself.
\end {proposition}

The proof of this statement is given in Section \ref {T_c1}. If a simplex $T$ is not an admissible domain for itself, we will need to narrow the class $W^2(T)$ to the set of functions that can be extended to some admissible domain $P$ for $T$ while remaining in the class $W^2(P)$.

\section {Main results}\label {S4results}

Throughout the paper we will denote by ${\bf a}\cdot {\bf b}$ the dot-product of vectors ${\bf a},{\bf b}\in \RR^d$ to distinguish it from the product of numbers $ab$, where the dot will be omitted.

We obtain the following optimality result.

\begin {theorem}\label {Th1}
Let $T$ be a non-degenerate simplex in $\RR^d$ with vertices ${\bf v}_0$, ${\bf v}_1$, $\ldots, {\bf v}_d$ and the circumcenter ${\bf c}$ and let $P$ be an admissible domain for $T$. Then the recovery method
\begin {equation}\label {rec1}
p_f({\bf w}):=\sum\limits_{i=0}^{d}b_i(\varphi({\bf w}))\(f({\bf v}_i)+\nabla f({\bf v}_i)\cdot \({\bf w}-\varphi({\bf w})+\frac 12 ({\bf c}-{\bf v}_i)\)\)
\end {equation}
is optimal in the sense of Problem \ref {P1} for recovery of $f({\bf w})$ on the class $W^2(P)$ for any ${\bf w}\in P$. That is, method \eqref {rec1} is pointwise optimal on $W^2(P)$. Moreover, for every point ${\bf w}\in P$,
\begin {equation}\label {h}
e_{{\bf w}}(W^2(P))=\mathcal E({\bf w}):=\frac {R^2}{4}+\left|{\bf w}-\varphi({\bf w})\right|^2-\frac {1}{2}\left|{\bf w}-{\bf c}\right|^2,
\end {equation}
where $R$ is the circumradius of the simplex $T$. Method $p_f$ is also optimal on $W^2(P)$ for global recovery in the sense of Problem \ref {P2} in any monotone seminorm $\|\cdot\|$ on $C(P)$ and 
\begin {equation}\label {h1}
R(W^2(P);\|\cdot\|)=\|\mathcal E\|.
\end {equation}
\end {theorem}
Notice that the recovery spline \eqref {rec1} is defined on all of $\RR^d$ even though the function $f$ is defined only on $P$. Thus, method \eqref {rec1} can also be used for extrapolation.

Theorem \ref {Th1} holds for $P=\RR^d$ no matter what simplex $T$ is given.
We also remark that the case $d=1$ of Theorem \ref {Th1} is the known result by Bojanov \cite {Boj1975} that we described above. 

Since $W^2(P)$ is a convex and symmetric class of functions,
by the well known Smolyak's lemma (see \cite {Smo65} or \cite {Bah71}), among optimal algorithms for recovery of $f({\bf w})$ on the class $W^2(P)$ there is a linear one; i.e., with a linear mapping $\Psi$ in \eqref {PSI}. Such is algorithm in \eqref {rec1}.

Theorem \ref {Th1} and Proposition \ref {T_c} immediately imply the following result.

\begin {corollary}\label{wcT}
Let $T$ be a non-degenerate completely well-centered simplex in $\RR^d$. Then the method $p_f$ given by \eqref {rec1} is optimal in the sense of Problem \ref {P1} for recovery of $f({\bf w})$ on the class $W^2(T)$ for any ${\bf w}\in T$ with $e_{{\bf w}}(W^2(T))=\mathcal E({\bf w})$. It is also optimal in the sense of Problem \ref {P2} for global recovery of the class $W^2(T)$ in any monotone seminorm $\|\cdot\|$ on $C(T)$ with $R(W^2(T);\|\cdot\|)=\|\mathcal E\|$.
\end {corollary}

A similar pointwise optimal recovery problem for general configurations of nodes is studied by the author in another paper, which is in preparation. Other optimal recovery problems such as recovery of a linear positive functional or a linear positive operator on the class of twice differentiable multivariate functions are also considered there.

An important ingredient of the proof of Theorem \ref {Th1} is the following statement describing the properties of the error function $\mathcal E$ which can be considered as the multivariate analogue of the Euler spline $\phi_2$ and may be useful for many extremal problems in Multivariate Approximation Theory.

\begin {theorem}\label {Th2}
Let $T$ be a non-degenerate simplex in $\RR^d$ with vertices ${\bf v}_0$, ${\bf v}_1$, $\ldots, {\bf v}_d$, circumcenter ${\bf c}$ and circumradius $R$. Then the function 
$$
\mathcal E({\bf x})=\frac {R^2}{4}+\left|{\bf x}-\varphi({\bf x})\right|^2-\frac {1}{2}\left|{\bf x}-{\bf c}\right|^2
$$
is $C^1$-continuous on $\RR^d$ with the gradient $\nabla \mathcal E({\bf x})={\bf x}+{\bf c}-2\varphi(\bf x)$, ${\bf x}\in \RR^d$. Furthermore, the gradient $\nabla \mathcal E$ satisfies the inequality
\begin {equation}\label {gradient}
\left|\nabla \mathcal E({\bf x})-\nabla \mathcal E({\bf y})\right|\leq \left|{\bf x}-{\bf y}\right|,\ \ \ {\bf x},{\bf y}\in \RR^d,
\end {equation}
and hence, $\mathcal E\in W^2(\RR^d)$. Finally, $\mathcal E({\bf v}_i)=0$ and $\nabla \mathcal E({\bf v}_i)={\bf 0}$, $i=0,1,\ldots,d$.
\end {theorem}
More properties of the function $\mathcal E$ are described in Theorems \ref {E} and \ref {nd} further in the text.

Since $W^2(P)$ is a convex and centrally symmetric class of functions and 
$$
\sup\limits_{f\in W^2(P)\atop I(f)={\bf 0}}f({\bf w})\leq \frac {\left|{\bf w}-{\bf v}_0\right|^2}{2}<\infty,\ \ \ {\bf w}\in P,
$$
by a standard argument from \cite {Smo65} (see also \cite {Bah71}), 
for every point ${\bf w}\in P$, we have
$$
e_{\bf w}(W^2(P))=\sup\limits_{f\in \mathcal W^2(P)\atop I(f)={\bf 0}}f({\bf w}).
$$
This fact and Theorem \ref {Th1} immediately imply the following result.
\begin {corollary}\label {q1'}
Let $T$ be a non-degenerate simplex in $\RR^d$ with vertices ${\bf v}_0$, ${\bf v}_1$, $\ldots, {\bf v}_d$, $P$ be an admissible domain for $T$, and ${\bf w}\in P$ be arbitrary point. Then
$$
\sup\limits_{f\in W^2(P)\atop I(f)={\bf 0}}f({\bf w})=\mathcal E({\bf w}). 
$$
\end {corollary}

\section {Further analysis of the optimal partition and of the error function}\label {S5p}

We devote this section to discussing some properties of the error function $\mathcal E({\bf x})$ and of the underlying partition. We say that two non-empty sets $A$ and $B$ in $\RR^d$ are mutually orthogonal if for every ${\bf a}_1,{\bf a}_2\in A$ and every ${\bf b}_1,{\bf b}_2\in B$, there holds $({\bf a}_1-{\bf a}_2)\cdot ({\bf b}_1-{\bf b}_2)=0$. Given two mutually orthogonal sets $A$ and $B$, we denote their direct sum by
$$
A\oplus B:=\left\{{\bf a}+{\bf b} : {\bf a}\in A,\ {\bf b}\in B \right\}
$$
(since the sets $A$ and $B$ are mutually orthogonal, every ${\bf x}\in A\oplus B$ has a unique representation ${\bf x}={\bf a}+{\bf b}$, where ${\bf a}\in A$ and ${\bf b}\in B$). Let also $\alpha A:=\{\alpha{\bf a} : {\bf a}\in A\}$, $\alpha>0$.

Given a set $X=\{{\bf v}_0,\ldots,{\bf v}_d\}\subset \RR^d$ in general position; i.e., not contained in any hyperplane,
let 
$$
W_i:=\{{\bf x}\in \RR^d : \left|{\bf x}-{\bf v}_i\right|\leq  \left|{\bf x}-{\bf v}_j\right|,\ j\in \Sigma\setminus \{i\}\}, \ \ \ i=0,1,\ldots,d,
$$
be the Voronoi cell of the point ${\bf v}_i$ with respect to the set $X$. For $I\subset \Sigma$, $I\neq \emptyset$, we denote 
$$
F_I:=\bigcap_{i\in I} W_i=\{{\bf x}\in \RR^d : \left|{\bf x}-{\bf v}_i\right|={\rm dist}({\bf x},X),\ i\in I\}.
$$ 
Notice that, for example, $F_\Sigma=\{{\bf c}\}$.

We start by showing that sets 
$$
\varphi^{-1}(G_I)=\{{\bf x}\in \RR^d : \varphi({\bf x})\in G_I\}, \ \ \ I\subset \Sigma, \ \ I\neq \emptyset, 
$$
form a partition of $\RR^d$ with the closure of each of them being a convex $d$-dimensional polyhedron (unbounded for $I\neq \Sigma$). Let 
$$
D_I:=\{{\bf x}\in T : \beta_i({\bf x})>0,\    i\in I, \ \ {\rm and}\ \ \beta_i({\bf x})=0,\ i\notin I \},
$$ 
where $\beta_i({\bf x})$, $i=0,1,\ldots,d$, are the barycentric coordinates of a point ${\bf x}$ with respect to the simplex~$T$.
\begin {proposition}\label {partition}
For every non-empty subset $I\subset \Sigma$, the sets $G_I$ and $F_I$ are mutually orthogonal, the set $\varphi^{-1}(G_I)$ is a convex set with non-empty interior, and
\begin {equation}\label {p2}
\varphi^{-1}(G_I)=G_I\oplus \frac 12(F_I-{\bf c})=\frac 12 (D_I\oplus F_I).
\end {equation}
\end {proposition}

Proposition \ref {partition} is proved in Section \ref {S10}.
It implies the following. Every set $A_{ij}:=\{{\bf x}\in \RR^d : \left|{\bf x}-{\bf v}_i\right|\leq \left|{\bf x}-{\bf v}_j\right|\}$, $i\neq j$, is a closed half-space in $\RR^d$ defined by the inequality ${\bf n}_{ij}\cdot {\bf x}\leq b_{ij}$, where ${\bf n}_{ij}=2({\bf v}_j-{\bf v}_i)$ and $b_{ij}={\bf v}_j^2-{\bf v}_i^2$. Since each set $F_I$ can be defined as the intersection of finitely many closed half-spaces $A_{ij}$, we have a matrix $B$ and a vector ${\bf b}$ such that $F_I=\{{\bf x}\in \RR^d : B{\bf x}\leq {\bf b}\}$. Since ${\bf n}_{ij}\cdot {\bf c}=b_{ij}$ for all $i\neq j$, we have $B{\bf c}={\bf b}$ and hence, $F_I-{\bf c}=\{{\bf x}\in \RR^d : B{\bf x}\leq {\bf 0}\}$. Thus, by definition, the set $F_I-{\bf c}$ is a polyhedral cone. Since $F_I-{\bf c}$ is a cone, we have $F_I-{\bf c}=(1/2)(F_I-{\bf c})$. In view of Proposition~\ref {partition}, we have $\OL {\varphi^{-1}(G_I)}=\OL G_I \oplus (1/2)(F_I-{\bf c})$, where $\OL {A}$ stands for the closure of the set $A$. Since $\OL G_I$ is a convex polytope and $(1/2)(F_I-{\bf c})$ is a polyhedral cone, the set $\OL {\varphi^{-1}(G_I)}$ is a convex polyherdron by the Minkowski-Weyl decomposition theorem for polyhedra (see e.g. \cite [Theorem 6.18]{CinVanRCA}). By Proposition \ref {partition}, each polyhedron $\OL {\varphi^{-1}(G_I)}$ has a non-empty interior. The union of the polyhedra $\OL {\varphi^{-1}(G_I)}$, $I\subset \Sigma$, $I\neq \emptyset$, is all of $\RR^d$ and, since the sets $G_I$ are pairwise disjoint, so are the interiors of the polyhedra $\OL {\varphi^{-1}(G_I)}$.

For a given subset $I\subset \Sigma$, denote by $[I,\Sigma]$ the collection of subsets $M$ of $\Sigma$ that contain the set $I$ (that is $I\subset M\subset \Sigma$). Let ${\bf u}_i^M:=\frac 12 ({\bf c}_M+{\bf v}_i)$, $i\in M$, $M\subset \Sigma$, $M\neq \emptyset$, where we recall that ${\bf c}_M$ is the circumcenter of the simplex $T_M$.
\begin {theorem}\label {well_c}
Let $T$ be a non-degenerate completely well-centered simplex in $\RR^d$. Then for every subset $I\subset \Sigma$, $I\neq \emptyset$, we have
$$
\varphi^{-1}(G_I)\cap T\!=\!\left\{\!\sum_{M\in [I,\Sigma]}\!\!\!\alpha_M\! \sum\limits_{i\in I}\beta_i {\bf u}_i^M : \sum_{i\in I}\!\beta_i =1; \!\!\!\sum_{M\in [I,\Sigma]}\!\!\!\alpha_M=1;\ \beta_i>0, i\in I; \alpha_M\geq 0,M\!\in\! [I,\Sigma]\right\}\!.
$$
\end {theorem}
Theorem \ref {well_c} follows immediately from Theorem \ref {set} proved in Section \ref {T_c1}.

We next show that the restriction of the function $\mathcal E({\bf x})$ to every convex polyhedron $\OL {\varphi^{-1}(G_I)}$ is a quadratic polynomial. Recall that for a given non-empty subset $I\subset \Sigma$, $Q_I$ denotes the affine subspace spanned by the set of vertices $\{{\bf v}_i\}_{i\in I}$. Recall that $H_I$ is the set of all points in $\RR^d$ that are equidistant from the vertices ${\bf v}_i$, $i\in I$. Recall also that ${\bf c}_I\in Q_I$ is the circumcenter of the simplex $T_I$ with vertices ${\bf v}_i$, $i\in I$, and that $R_I$ is the circumradius of $T_I$. Notice that $m:=\# I-1$ is the dimension of $Q_I$ and that $H_I$ is a $(d-m)$-dimensional affine subspace which passes through the point ${\bf c}_I$ and is mutually orthogonal with $Q_I$ (see Lemma \ref {L10a}).

\begin {theorem}\label {E}
Let $T$ be a non-degenerate simplex in $\RR^d$ and let $I\subset \Sigma$ be a non-empty set. Then
\begin {equation}\label {E1}
\mathcal E({\bf x})=\frac {R_I^2}{4}-\frac 12 \({\rm dist}({\bf x},H_I)\)^2+\frac 12 \({\rm dist}({\bf x},Q_I)\)^2, \ \ \ {\bf x}\in \OL {\varphi^{-1}(G_I)}.
\end {equation}
Furthermore, the restriction of the gradient $\nabla \mathcal E$ of the function $\mathcal E$ to the set $\OL {\varphi^{-1}(G_I)}$ is an isometry; i.e., 
\begin {equation}\label {E2}
\left|\nabla \mathcal E({\bf x})-\nabla \mathcal E({\bf y})\right|=\left|{\bf x}-{\bf y}\right|,\ \ \ {\bf x},{\bf y}\in \OL {\varphi^{-1}(G_I)}.
\end {equation}
\end {theorem}
Theorem \ref {E} is proved in Section \ref {S10}.

{\bf Remark.} If $L$ is any affine subspace of $\RR^d$ (of dimension $0\leq l<d$), then the function $g({\bf x}):=({\rm dist}({\bf x},L))^2$ is a non-negative quadratic polynomial. If, in particular, $l=0$, then $L$ consists of one point ${\bf z}_0$ and $g({\bf x})=({\bf x}-{\bf z}_0)^2$, ${\bf x}\in \RR^d$. If $l=d$, then $L=\RR^d$ and we simply have $g({\bf x})=0$.

{\bf Remark.} In the case $I=\Sigma$ in Theorem \ref {E}, we have $Q_I=\RR^d$, $H_I=\{{\bf c}\}$, $G_I={\rm int}\ \! U$, $\varphi^{-1}(G_I)={\rm int}\ \! U$, and $\mathcal E({\bf x})=\frac {R^2}{4}-\frac 12 \left|{\bf x}-{\bf c}\right|^2$, $x\in U$. If $I=\{i\}$, then $Q_I=\{{\bf v}_i\}$, $R_I=0$, $H_I=\RR^d$, $G_I=\{{\bf u}_i\}$, $\varphi^{-1}(G_I)=\{{\bf u}_i\}\oplus (1/2)(W_i-{\bf c})=({\bf v}_i+W_i)/2$ (see Proposition \ref {partition}). Consequently, $\mathcal E({\bf x})=\frac 12 ({\bf x}-{\bf v}_i)^2$, ${\bf x}\in ({\bf v}_i+W_i)/2$. In all other cases $\mathcal E({\bf x})$ is a difference of two non-negative quadratic polynomials.

We next give sufficient conditions under which the normal derivative of the error function $\mathcal E$ vanishes on the interior of every facet (that is, $(d-1)$-dimensional face) of $T$.
\begin {theorem}\label {nd}
Let $T$ be a non-degenerate simplex in $\RR^d$ that contains its circumcenter and is such that every facet of $T$ contains the orthogonal projection of the opposite vertex of $T$.
If ${\bf x}$ is any point in the relative interior of any facet of $T$ and ${\bf n}$ is a unit normal vector to that facet, then $\frac {\partial \mathcal E}{\partial {\bf n}}({\bf x})=0$.
\end {theorem}

\section {Auxiliary statements}\label {aux}

We start with the following well-known statement. For completeness, we present its proof in the Appendix.
\begin {lemma}\label {distconv}
Let $K$ be a non-empty compact convex set in $\RR^d$ and ${\bf x}\in \RR^d$ be arbitrary point. Then there is a unique point ${\bf y}$ in $K$ closest to ${\bf x}$. Furthermore, the inequality
\begin {equation}\label {q}
({\bf x}-{\bf y})~\cdot~ ({\bf z}-{\bf y})\leq 0,\ \ \ {\bf z}\in K,
\end{equation}
holds if and only if ${\bf y}$ is the point in $K$ closest to ${\bf x}$.
\end {lemma}
We next establish the Lipschitz continuity of the function $\varphi$.
\begin {lemma}\label {phiLip}
For every ${\bf x},{\bf y}\in \RR^d$, there holds $\left|\varphi({\bf x})-\varphi({\bf y})\right|\leq \left|{\bf x}-{\bf y}\right|$.
\end {lemma}
\begin {proof}
Since $\varphi$ maps every point in $\RR^d$ into the closest point in $U$, by Lemma \ref {distconv} for any ${\bf x},{\bf y}\in \RR^d$, we will have
\begin {equation}\label {x_phi}
({\bf x}-\varphi({\bf x}))\cdot (\varphi({\bf y})-\varphi({\bf x}))\leq 0\ \ \ {\rm and}\ \ \ ({\bf y}-\varphi({\bf y}))\cdot (\varphi({\bf x})-\varphi({\bf y}))\leq 0.
\end {equation}
Using also the inequality
$$
({\bf x}-{\bf y})\cdot (\varphi({\bf x})-\varphi({\bf y}))\leq \frac 12 \((\varphi({\bf x})-\varphi({\bf y}))^2+({\bf x}-{\bf y})^2\)
$$
we will obtain
$$
({\bf x}-{\bf y})^2-(\varphi({\bf x})-\varphi({\bf y}))^2=({\bf x}-{\bf y})^2+({\bf x}-\varphi({\bf x}))\cdot (\varphi({\bf x})-\varphi({\bf y}))
$$
$$
-({\bf y}-\varphi({\bf y}))\cdot (\varphi({\bf x})-\varphi({\bf y}))-({\bf x}-{\bf y})\cdot (\varphi({\bf x})-\varphi({\bf y}))
$$
$$
\geq ({\bf x}-{\bf y})^2-({\bf x}-{\bf y})\cdot (\varphi({\bf x})-\varphi({\bf y}))
$$
$$
\geq \frac 12 \(({\bf x}-{\bf y})^2-(\varphi({\bf x})-\varphi({\bf y}))^2\).
$$
Consequently, $({\bf x}-{\bf y})^2-(\varphi({\bf x})-\varphi({\bf y}))^2\geq 0$ and the assertion of the lemma follows.
\end {proof}
We will also need the following statement.
\begin {proposition}\label {W2}
Let $K$ be a convex set in $\RR^d$ with non-empty interior. 
If $f\in C^1(K)$ and $\left|\nabla f({\bf x})-\nabla f({\bf y})\right|\leq \left|{\bf x}-{\bf y}\right|$, ${\bf x},{\bf y}\in K$, then $f\in W^2(K)$.
\end {proposition}
\begin {proof}
Let $g(t)=f\({\bf x}+\frac tr ({\bf y}-{\bf x})\)$, $t\in [0,r]$, where $r=\left|{\bf x}-{\bf y}\right|$ and ${\bf x},{\bf y}\in K$ are arbitrary distinct points. Then $g'(t)=\nabla f \({\bf x}+\frac {t}{r}({\bf y}-{\bf x})\)\cdot \frac {{\bf y}-{\bf x}}{r}$ and for every $u,v\in [0,r]$, we have
$$
\left|g'(u)-g'(v)\right|=\left|\nabla f\({\bf x}+\frac {u} {r} ({\bf y}-{\bf x})\)\cdot \frac {{\bf y}-{\bf x}}{r}-\nabla f\({\bf x}+\frac v r ({\bf y}-{\bf x})\)\cdot \frac {{\bf y}-{\bf x}}{r}\right|
$$
$$
\leq \left|\nabla f\({\bf x}+\frac {u} {r} ({\bf y}-{\bf x})\)-\nabla f\({\bf x}+\frac v r ({\bf y}-{\bf x})\)\right|\left|\frac {{\bf y}-{\bf x}}{r}\right|
$$
$$
\leq \left|\frac {u}{r}({\bf y}-{\bf x})-\frac {v}{r}({\bf y}-{\bf x})\right|=\left|u-v\right|,
$$
which implies that $g\in W^2_\infty[0,r]$. Consequently, $g\in W^2(K)$.
\end {proof}

\section {Proof of Proposition \ref {T_c} and of Theorem \ref {well_c}}\label {T_c1}

To prove Theorem \ref {well_c} and Proposition \ref {T_c}, we will first establish the following auxiliary statements. Some of the proofs in this section are straightforward but we keep them for completeness.
\begin {lemma}\label {T_J}
Let $T_J$ be a simplex in $\RR^d$ with an affinely independent set of vertices $\{{\bf w}_i\}_{i\in J}$ ($2\leq \# J\leq d+1$). Let also $j\in J$ be any index and let $I:=J\setminus \{j\}$. If $T_J$ contains its circumcenter, then
$$
R_I\leq \left|{\bf c}_I-{\bf w}_j\right|,
$$
where ${\bf c}_I$ denotes the circumcenter and $R_I$ denotes the circumradius of the face $T_I$ of $T_J$.
\end {lemma}
\begin {proof}
Assume to the contrary that $R_I>\left|{\bf c}_I-{\bf w}_j\right|$. Fix some index $m\in I$. Let ${\bf n}$ be a non-zero vector in the subspace $H_1:={\rm span}\{{\bf w}_i-{\bf w}_m\}_{i\in J\setminus \{m\}}$ that is orthogonal to every vector ${\bf w}_i-{\bf w}_m$, $i\in I\setminus \{m\}$. Then, clearly, ${\bf n}\cdot ({\bf w}_j-{\bf w}_m)\neq 0$ and we will choose the direction of ${\bf n}$ so that ${\bf n}\cdot ({\bf w}_j-{\bf w}_m)< 0$. Notice that for every $i\in I$,
$$
{\bf n}\cdot ({\bf w}_i-{\bf w}_j)={\bf n}\cdot\( ({\bf w}_i-{\bf w}_m)+({\bf w}_m-{\bf w}_j)\)={\bf n}\cdot ({\bf w}_m-{\bf w}_j).
$$
The circumcenter ${\bf c}_I$ can be written as ${\bf c}_I=\sum_{i\in I}\alpha_i{\bf w}_i$, where $\sum_{i\in I}\alpha_i=1$. Then
$$
{\bf n}\cdot ({\bf c}_I-{\bf w}_j)={\bf n}\cdot\(\sum_{i\in I}\alpha_i({\bf w}_i-{\bf w}_j)\)=\sum_{i\in I}\alpha_i {\bf n}\cdot ({\bf w}_i-{\bf w}_j)
$$
$$
=\sum_{i\in I}\alpha_i {\bf n}\cdot ({\bf w}_m-{\bf w}_j)= {\bf n}\cdot ({\bf w}_m-{\bf w}_j)>0.
$$
It is also clear that
$$
{\bf n}\cdot ({\bf c}_I-{\bf w}_i)={\bf n}\cdot ({\bf c}_I-{\bf w}_j)+{\bf n}\cdot ({\bf w}_j-{\bf w}_i)={\bf n}\cdot ({\bf w}_m-{\bf w}_j)+{\bf n}\cdot ({\bf w}_j-{\bf w}_m)=0,\ \ \ i\in I.
$$
Letting ${\bf x}_t:={\bf c}_I+t{\bf n}$, $t\geq 0$ (notice that ${\bf x}_t\in Q_J$), we have
$$
\left|{\bf x}_t-{\bf w}_i\right|^2=t^2{\bf n}^2+2t{\bf n}\cdot ({\bf c}_I-{\bf w}_i)+({\bf c}_I-{\bf w}_i)^2=t^2{\bf n}^2+R_I^2,\ \ \ i\in I.
$$
Then the distance $\left|{\bf x}_t-{\bf w}_i\right|$ is independent on $i\in I$. Notice that
$$
g(t):=\left|{\bf x}_t-{\bf w}_j\right|^2-\left|{\bf x}_t-{\bf w}_m\right|^2=2t{\bf n}\cdot ({\bf c}_I-{\bf w}_j)+({\bf c}_I-{\bf w}_j)^2-R_I^2.
$$
By the contrary assumption, $g(0)<0$. Furthermore, $\lim_{t\to+\infty} g(t)=+\infty$. Hence, there exists $t_1>0$ such that $g(t_1)=0$. Then 
$$
\left|{\bf x}_{t_1}-{\bf w}_j\right|=\left|{\bf x}_{t_1}-{\bf w}_m\right|=t_1^2{\bf n}^2+R_I^2=\left|{\bf x}_{t_1}-{\bf w}_i\right|, \ \ \ i\in I. 
$$
Hence, the point ${\bf x}_{t_1}\in Q_J$ coincides with the circumcenter ${\bf c}_J$ of the simplex $T_J$. 

Let $\{\beta_i\}_{i\in J}$ be the barycentric coordinates of the point ${\bf c}_J$ relative to the simplex $T_J$. By assumption, $\beta_j\geq 0$. Then by the choice of the vector ${\bf n}$, we have
$$
{\bf n}\cdot ({\bf c}_J-{\bf w}_m)=\sum_{i\in J}\beta_i {\bf n}\cdot ({\bf w}_i-{\bf w}_m)=\beta_j{\bf n}\cdot ({\bf w}_j-{\bf w}_m)\leq 0
$$
while the fact that ${\bf c}_J={\bf x}_{t_1}$ implies that
$$
{\bf n}\cdot ({\bf c}_J-{\bf w}_m)={\bf n}\cdot ({\bf x}_{t_1}-{\bf w}_m)={\bf n}\cdot ({\bf c}_I-{\bf w}_m)+t_1{\bf n}^2=t_1{\bf n}^2>0.
$$
This contradiction shows that $R_I\leq \left|{\bf c}_I-{\bf w}_j\right|$. 
\end {proof}

\begin {corollary}\label {dot_p}
Let $T_L$ be a simplex in $\RR^d$ with an affinely independent set of vertices $\{{\bf w}_i\}_{i\in L}$ and let $M$ be a non-empty proper subset of $L$. Assume that every face $T_J$ of $T_L$ such that $J\supset M$ and $\# J=\# M+1$, contains its circumcenter. Then for every $i\in M$ and $j\in L$, 
$$
({\bf c}_M-{\bf c}_L)\cdot ({\bf w}_j-{\bf w}_i)\leq 0.
$$
\end {corollary}
\begin {proof}
Let $i\in M$ be arbitrary. If $j\in M$, we have $({\bf c}_M-{\bf w}_i)^2=({\bf c}_M-{\bf w}_j)^2$. If $j\in L\setminus M$, we let $J_1:=M\cup \{j\}$. By assumption, the simplex $T_{J_1}$ contains its circumcenter. In view of Lemma \ref {T_J}, we have $({\bf c}_M-{\bf w}_i)^2=R_M^2\leq ({\bf c}_M-{\bf w}_j)^2$. Thus, for every $i\in M$ and $j\in L$, we have the following set of inequalities (since $i,j\in L$, we use the equality $\left|{\bf c}_L-{\bf w}_i\right|=\left|{\bf c}_L-{\bf w}_j\right|$)
$$
({\bf c}_M-{\bf w}_i)^2\leq ({\bf c}_M-{\bf w}_j)^2
$$
$$
({\bf c}_M-{\bf c}_L+{\bf c}_L-{\bf w}_i)^2\leq ({\bf c}_M-{\bf c}_L+{\bf c}_L-{\bf w}_j)^2;
$$
$$
({\bf c}_M-{\bf c}_L)^2+2({\bf c}_M-{\bf c}_L)\cdot ({\bf c}_L-{\bf w}_i)+({\bf c}_L-{\bf w}_i)^2\leq ({\bf c}_M-{\bf c}_L)^2+2({\bf c}_M-{\bf c}_L)\cdot ({\bf c}_L-{\bf w}_j)+({\bf c}_L-{\bf w}_j)^2;
$$
$$
({\bf c}_M-{\bf c}_L)\cdot ({\bf c}_L-{\bf w}_i)\leq ({\bf c}_M-{\bf c}_L)\cdot ({\bf c}_L-{\bf w}_j);
$$
$$
({\bf c}_M-{\bf c}_L)\cdot ({\bf w}_j-{\bf w}_i)\leq 0,
$$
which implies the assertion of the corollary. 
\end {proof}

Before proving the next statement we introduce some auxiliary notation. Let $L$ be a set of indices, $\# L\leq d+1$, and let $T_L$ be a simplex in $\RR^d$ with an affinely independent set of vertices $\{{\bf w}_i\}_{i\in L}$. Denote 
$$
{\bf u}^L_i:=\frac 12 ({\bf c}_L+{\bf w}_i), \ \ i\in L,\ \ \ {\rm and \ let}\ \ \  U_L:=\frac 12 ({\bf c}_L+T_L).
$$
Then $U_L$ is the simplex with the set of vertices $\{{\bf u}^L_i\}_{i\in L}$.  

Recall that $Q_L$ denotes the affine subspace spanned by the point set $\{{\bf w}_i\}_{i\in L}$; i.e.,  
$$
Q_L=\left\{\sum_{i\in L}c_i{\bf w}_i : \sum_{i\in L}c_i=1\right\}.
$$
Let also $\varphi_L:Q_L\to U_L$ be the mapping such that for every ${\bf x}\in Q_L$, $\varphi_L({\bf x})$ is the point in $U_L$ closest to ${\bf x}$.

\begin {lemma}\label {conv_phi}
Let $T_L$ be a simplex in $\RR^d$ with affinely independent set of vertices and let ${\bf a}\in U_L$ be arbitrary point. Then the set $\varphi_L^{-1}(\{{\bf a}\})=\{{\bf x}\in Q_L : \varphi_L({\bf x})={\bf a}\}$ is convex.
\end {lemma}
\begin {proof}
Let ${\bf x},{\bf y}\in \varphi_L^{-1}(\{{\bf a}\})$ be arbitrary points and let $t\in [0,1]$.
Since ${\bf a}$ is the point in $U_L$ closest to ${\bf x}$ and to ${\bf y}$, by Lemma \ref {distconv}, for every ${\bf u}\in U_L$, we have $({\bf x}-{\bf a})\cdot ({\bf u}-{\bf a})\leq 0$ and $({\bf y}-{\bf a})\cdot ({\bf u}-{\bf a})\leq 0$. Then
$$
(t{\bf x}+(1-t){\bf y}-{\bf a})\cdot ({\bf u}-{\bf a})=t({\bf x}-{\bf a})\cdot ({\bf u}-{\bf a})+(1-t)({\bf y}-{\bf a})\cdot ({\bf u}-{\bf a})\leq 0.
$$
Applying again Lemma \ref {distconv}, we obtain that $\varphi_L(t{\bf x}+(1-t){\bf y})={\bf a}$. Then $t{\bf x}+(1-t){\bf y}\in \varphi_L^{-1}(\{{\bf a}\})$ and the convexity of the set $\varphi_L^{-1}(\{{\bf a}\})$ follows.
\end {proof}

\medskip

Recall that $[I,L]$, $I\subset L$, is the collection of all subsets of $L$ that contain the set $I$. 

\begin {lemma}\label {phi_J}
Let $L$ be a set of indices, $2\leq\# L\leq d+1$, and let $T_L$ be a simplex in $\RR^d$ with an affinely independent set of vertices $\{{\bf w}_i\}_{i\in L}$. Let $M$ be a non-empty proper subset of $L$ such that every simplex $T_J$, $J\in [M,L]$, $J\neq M$, contains its circumcenter. Then for every point ${\bf z}:=\sum_{i\in M}\beta_i{\bf u}^M_i$ such that $\beta_i\geq 0$, $i\in M$, and $\sum_{i\in M}\beta_i=1$, there holds $\varphi_L({\bf z})=\sum_{i\in M}\beta_i{\bf u}^L_i$.
\end {lemma}
\begin {proof}
Let ${\bf y}:=\sum_{i\in M}\beta_i{\bf u}^L_i$.
Notice that ${\bf u}^M_i-{\bf u}_i^L=\frac 12 ({\bf c}_M+{\bf w}_i)-\frac 12 ({\bf c}_L+{\bf w}_i)=\frac 12 ({\bf c}_M-{\bf c}_L)$, $i\in M$. For every $i,l\in L$, we also have ${\bf u}^L_l-{\bf u}^L_i=\frac 12 ({\bf w}_l-{\bf w}_i)$. By Corollary \ref {dot_p}, for every $l\in L$, we obtain
$$
\({\bf z}-{\bf y}\)\cdot \({\bf u}_l^L-{\bf y}\)=\(\sum_{i\in M}\beta_i({\bf u}_i^M-{\bf u}_i^L)\)\cdot \(\sum_{i\in M}\beta_i({\bf u}_l^L-{\bf u}_i^L)\)
$$
$$
=\frac 14 ({\bf c}_M-{\bf c}_L)\cdot \(\sum_{i\in M}\beta_i({\bf w}_l-{\bf w}_i)\)=\frac 14\sum_{i\in M}\beta_i({\bf c}_M-{\bf c}_L)\cdot ({\bf w}_l-{\bf w}_i)\leq 0.
$$
Since every point ${\bf u}\in U_L$ can be written as a convex combination of the points ${\bf u}^L_l$, $l\in L$, we have $({\bf z}-{\bf y})\cdot ({\bf u}-{\bf y})\leq 0$, ${\bf u}\in U_L$. Then by Lemma \ref {distconv}, we have $\varphi_L({\bf z})={\bf y}$. 
\end {proof}

For every non-empty subset $I\subset L$, denote
$$
G^L_I:=\left\{\sum_{i\in L}\beta_i {\bf u}^L_i : \sum _{i\in L} \beta_i =1; \  \beta_i>0,\ i\in I;\ \beta_i=0,\ i\in L\setminus I\right\}
$$
and let
$$
V^L_I:=\left\{\sum_{M\in [I,L]}\!\!\alpha_M\sum_{i\in I}\beta_i{\bf u}_i^M : \sum_{i\in I}\beta_i=1;\ \beta_i>0,\ i\in I;\!\! \sum_{M\in [I,L]}\!\!\alpha_M=1;\ \alpha_M\geq 0,\ M\in [I,L]\right\}.
$$

\begin {lemma}\label {subset}
Let $L$ be a set of indices, $2\leq \# L\leq d+1$, and let $T_L$ be a simplex in $\RR^d$ with an affinely independent set of vertices $\{{\bf w}_i\}_{i\in L}$. Let $I$ be a non-empty proper subset of $L$. If for every set $M\in [I,L]$, $M\neq I$, the simplex $T_M$ contains its circumcenter, then $V^L_I\subset \varphi_L^{-1}(G^L_I)$. If, in addition, $T_I$ contains its circumcenter, then $V^L_I\subset \varphi_L^{-1}(G^L_I)\cap T_L$.
\end {lemma}
\begin {proof}
Every point ${\bf x}\in V^L_I$ has a representation ${\bf x}=\sum\limits_{M\in [I,L]}\!\!\alpha_M\sum\limits_{i\in I}\beta_i{\bf u}_i^M$, where $\sum_{i\in I}\beta_i=1$, $\beta_i>0$, $i\in I$, $\sum\limits_{M\in [I,L]}\!\!\alpha_M=1$, and $\alpha_M\geq 0$, $M\in [I,L]$. Then ${\bf x}\in Q_L$.

By Lemma \ref {phi_J}, for every $M\in [I,L]$, $M\neq L$, and the point ${\bf z}_M:=\sum_{i\in I}\beta_i{\bf u}^M_i\in Q_L$, we have $\varphi_L({\bf z}_M)=\sum_{i\in I}\beta_i{\bf u}^L_i$. Moreover, for  ${\bf z}_L:=\sum_{i\in I}\beta_i{\bf u}^L_i\in U_L$, we have $\varphi_L({\bf z}_L)={\bf z}_L$. Thus, ${\bf z}_M\in \varphi_L^{-1}(\sum_{i\in I}\beta_i{\bf u}_i^L)$ for any $M\in [I,L]$. Then by Lemma \ref {conv_phi}, we have
$$
\varphi_L({\bf x})=\varphi_L\(\sum\limits_{M\in [I,L]}\!\!\alpha_M{\bf z}_M\)=\sum_{i\in I}\beta_i{\bf u}^L_i\in G^L_I;
$$
i.e., ${\bf x}\in \varphi_L^{-1}(G^L_I)$.
Consequently, $V^L_I\subset \varphi_L^{-1}(G^L_I)$. 

If we assume that $T_I$ contains its circumcenter, then ${\bf c}_M\in T_M$, $M\in [I,L]$ and, hence, ${\bf u}^M_i=\frac 12 ({\bf c}_M+{\bf w}_i)\in T_M$, $i\in M$, $M\in [I,L]$. Then ${\bf z}_M\in T_M\subset T_L$, $M\in [I,L]$, and, consequently, ${\bf x}\in T_L$; i.e., $V^L_I\subset \varphi_L^{-1}(G^L_I)\cap T_L$. 
\end {proof}

\begin {theorem}\label {set}
Let $L$ be a set of indices, $2\leq \# L\leq d+1$, and let $T_L$ be a completely well-centered simplex in $\RR^d$ with an affinely independent set of vertices $\{{\bf w}_i\}_{i\in L}$. Then for every non-empty subset $I$ of $L$, there holds $V^L_I=\varphi_L^{-1}(G^L_I)\cap T_L$.
\end {theorem}
\begin {proof}
We will prove this result using induction on the cardinality of $L$. If $\#L=2$, without loss of generality, we can let $L=\{0,1\}$. For points ${\bf x},{\bf y}\in \RR^d$, let $[{\bf x},{\bf y}]=\{t{\bf x}+(1-t){\bf y} : t\in [0,1]\}$. Then ${\bf c}_L=\frac 12 ({\bf w}_0+{\bf w}_1)$ and $\varphi_L({\bf x})={\bf u}^L_i$, ${\bf x}\in [{\bf w}_i,{\bf u}_i^L]$, $i=0,1$. Consequently, 
$$
V^L_{\{i\}}=[{\bf w}_i,{\bf u}_i^L]=\varphi_L^{-1}(\{{\bf u}^L_i\})\cap [{\bf w}_0,{\bf w}_1]=\varphi_L^{-1}(G^L_{\{i\}})\cap T_L, \ \ \ i=0,1. 
$$
Moreover, 
$$
V_L^L=\{t{\bf u}^L_0+(1-t){\bf u}^L_1 : t\in (0,1)\}=G_L^L=\varphi_L^{-1}(G^L_L)=\varphi_L^{-1}(G^L_L)\cap T_L.
$$

Let now $2<k\leq d+1$ and assume that the assertion of the theorem holds for any index set $L$ with $2\leq \#L\leq k-1$. Choose any index set $L$ of cardinality $k$. Show that the theorem holds for $L$. Clearly, $V_L^L=G_L^L=\varphi_L^{-1}(G_L^L)\cap T_L$ and in view of Lemma \ref {subset}, we have $V_I^L\subset \varphi_L^{-1}(G_I^L)\cap T_L$ for every non-empty proper subset $I\subset L$. It remains to prove that $\varphi_L^{-1}(G^L_I)\cap T_L\subset V^L_I$ for every non-empty proper subset $I\subset L$.  

Choose any ${\bf x}\in \varphi_L^{-1}(G^L_I)\cap T_L$. If ${\bf x}\in G^L_I$, then, clearly, ${\bf x}\in V^L_I$. Assume that ${\bf x}\notin G^L_I$. Since $\varphi_L({\bf x})\in G^L_I$, we have $\varphi_L({\bf x})\neq {\bf x}$. Denote ${\bf z}_t:=\varphi_L({\bf x})+t({\bf x}-\varphi_L({\bf x}))$, $t\geq 1$, and let ${\bf x}=\sum_{i\in L}\gamma_i {\bf w}_i$ and $\varphi_L({\bf x})=\sum_{i\in L}\delta_i{\bf w}_i$ be the barycentric representations of the points ${\bf x}$ and $\varphi_L({\bf x})$ relative to the simplex $T_L$. Then ${\bf z}_t=\sum_{i\in L}(\delta_i+t(\gamma_i-\delta_i)){\bf w}_i$. 
Since $\varphi_L({\bf x})\neq {\bf x}$, we have $\gamma_i\neq \delta_i$ for some $i\in L$. Since $\sum_{i\in L}\gamma_i=\sum_{i\in L}\delta_i(=1)$, for some $v,m\in L$, we have $\gamma_v>\delta_v$ and $\gamma_m<\delta_m$.
Then the coordinate $\delta_m+t(\gamma_m-\delta_m)$ strictly decreases as $t$ increases. Let $t_0$ be the maximal value of $t$ in the interval $[1,\infty)$ such that $\min_{i\in L}\{\delta_i+t(\gamma_i-\delta_i)\}\geq 0$. Then $\delta_l+t_0(\gamma_l-\delta_l)=0$ for some $l\in L$. Let $J:=\{i\in L : \delta_i+t_0(\gamma_i-\delta_i)>0\}$. Since $\sum_{i\in L}(\delta_i+t_0(\gamma_i-\delta_i))=1$, we have $J\neq \emptyset$. Then $1\leq \# J\leq k-1$ and ${\bf z}_{t_0}\in T_J$.

In the case $\# J\geq 2$, let a non-empty subset $I_0\subset J$ be such that $\varphi_J({\bf z}_{t_0})\in G^J_{I_0}$. Since $T_J$ is also completely well-centered, by the induction assumption, we have $V^J_{I_0}=\varphi_J^{-1}(G^J_{I_0})\cap T_J$. Consequently, ${\bf z}_{t_0}\in V^J_{I_0}$. 
Then 
\begin {equation}\label {1'w}
{\bf z}_{t_0}=\sum\limits_{M\in [I_0,J]}\!\!\alpha_M\sum_{i\in I_0}\beta_i{\bf u}_i^M,
\end {equation}
where $\sum_{i\in I_0}\beta_i=1$, $\beta_i>0$, $i\in I_0$, $\sum\limits_{M\in [I_0,J]}\!\!\alpha_M=1$, and $\alpha_M\geq 0$, $M\in [I_0,J]$. 
By Lemma \ref {phi_J}, we have $\varphi_L\(\sum_{i\in I_0}\beta_i{\bf u}_i^M\)=\sum_{i\in I_0}\beta_i{\bf u}_i^L$, $M\in [I_0,J]$. Then by Lemma \ref {conv_phi}, we have $\varphi_L\({\bf z}_{t_0}\)=\sum_{i\in I_0}\beta_i{\bf u}_i^L\in G^L_{I_0}$. 

In the case $\# J=1$, since $\sum_{j\in L}(\delta_j+t_0(\gamma_j-\delta_j))=1$, we have $\delta_u+t_0(\gamma_u-\delta_u)=1$ for some $u\in L$ and hence, $J=\{u\}$ and ${\bf z}_{t_0}={\bf w}_u={\bf u}_u^{\{u\}}$. By Lemma \ref {phi_J}, we then have $\varphi_L({\bf z}_{t_0})={\bf u}_u^L\in G^L_{\{u\}}$. With $I_0=\{u\}$, $\beta_u=1$, and $\alpha_{\{u\}}=1$ we also obtain that $\varphi_L({\bf z}_{t_0})\in G^L_{I_0}$ and \eqref {1'w} holds.

In both cases, for every ${\bf u}\in U_L$, by Lemma \ref {distconv}, we have
$$
({\bf z}_{t_0}-\varphi_L({\bf x}))\cdot ({\bf u}-\varphi_L({\bf x}))=t_0({\bf x}-\varphi_L({\bf x}))\cdot ({\bf u}-\varphi_L({\bf x}))\leq 0.
$$ 
Hence, by Lemma \ref {distconv}, $\varphi_L({\bf z}_{t_0})=\varphi_L({\bf x})\in G^L_I$. Since the set $\{{\bf w}_i\}_{i\in L}$ is affinely independent, so is the set $\{{\bf u}_i^L\}_{i\in L}$. Since $\varphi_L({\bf z}_{t_0})\in G^L_{I_0}$, we have $I_0=I$. Thus, $\varphi_L({\bf x})=\sum_{i\in I}\beta_i{\bf u}_i^L$ and
$$
{\bf z}_{t_0}=\sum\limits_{M\in [I,J]}\!\!\alpha_M\sum_{i\in I}\beta_i{\bf u}_i^M.
$$
Since $1/t_0\in (0,1]$, we have
$$
{\bf x}=\frac {1}{t_0}{\bf z}_{t_0}+\(1-\frac {1}{t_0}\)\varphi_L({\bf x})=\sum\limits_{M\in [I,J]}\!\!\frac {\alpha_M}{t_0}\sum_{i\in I}\beta_i{\bf u}_i^M+\(1-\frac {1}{t_0}\)\sum_{i\in I}\beta_i{\bf u}_i^L\in V^L_I.
$$
Consequently, we obtain that $\varphi_L^{-1}(G^L_I)\cap T_L\subset V^L_I$. Then $\varphi_L^{-1}(G^L_I)\cap T_L=V^L_I$ for every $I\subset L$, $I\neq \emptyset$. Now the assertion of the theorem follows for any index set $L$ with $\#L=k$.
\end {proof}

\medskip

{\bf Proof of Theorem \ref {well_c}.} Theorem \ref {set} implies immediately Theorem \ref {well_c}.\hfill $\square$

\medskip

{\bf Proof of Proposition \ref {T_c}.} Let ${\bf x}\in T$ be arbitrary point and let $I\subset \Sigma$ be the non-empty subset such that $\varphi({\bf x})\in G_I$ (recall that we omit the index $\Sigma$ in the notation of $G^\Sigma_I$, $T_{\Sigma}$, and $\varphi_\Sigma$). We need to show that ${\bf x}-\varphi({\bf x})+{\bf u}_l\in T$, $l\in I$. By Theorem \ref {set}, we have $\varphi^{-1}(G_I)\cap T=V_I^\Sigma$. Then ${\bf x}\in V^\Sigma_I$ and hence, we have a representation 
$$
{\bf x}=\sum\limits_{M\in [I,\Sigma]}\!\!\alpha_M\sum_{i\in I}\beta_i{\bf u}_i^M,
$$
where ${\bf u}^M_i=\frac {1}{2}({\bf c}_M+{\bf v}_i)$, $M\in [I,\Sigma]$, $i\in I$, $\sum_{i\in I}\beta_i=1$, $\beta_i>0$, $i\in I$, $\sum\limits_{M\in [I,\Sigma]}\!\!\alpha_M=1$, and $\alpha_M\geq 0$, $M\in [I,\Sigma]$. By Lemma \ref {phi_J}, we have $\varphi\(\sum_{i\in I}\beta_i{\bf u}_i^M\)=\sum_{i\in I}\beta_i{\bf u}_i$, $M\in [I,\Sigma]$ (when $M=\Sigma$ this equality is trivial), and by Lemma \ref {conv_phi}, $\varphi({\bf x})=\sum_{i\in I}\beta_i{\bf u}_i$. Then
$$
{\bf x}-\varphi({\bf x})+{\bf u}_l=\sum\limits_{M\in [I,\Sigma]}\!\!\alpha_M\sum_{i\in I}\beta_i{\bf u}_i^M+{\bf u}_l-\sum_{i\in I}\beta_i{\bf u}_i=\sum\limits_{M\in [I,\Sigma]}\!\!\alpha_M\sum_{i\in I}\beta_i\({\bf u}_i^M+{\bf u}_l-{\bf u}_i\)
$$
$$
=\frac 12\sum\limits_{M\in [I,\Sigma]}\!\!\alpha_M\sum_{i\in I}\beta_i\({\bf c}_M+{\bf v}_l\)=\sum\limits_{M\in [I,\Sigma]}\!\!\alpha_M{\bf u}_l^M.
$$
By assumption, ${\bf c}_M\in T$, $M\in [I,\Sigma]$. Consequently, ${\bf u}^M_l\in T$, $M\in [I,\Sigma]$, and we obtain that ${\bf x}-\varphi({\bf x})+{\bf u}_l\in T$, $l\in I$. By Definition \ref {D3.1}, the simplex $T$ is an admissible domain for itself. \hfill $\square$

\section {Proof of Theorem \ref {Th2}}\label {S8_4.3}

It will be more convenient to prove Theorem \ref {Th2} first.
We start by showing that 
\begin {equation}\label {grad_e}
\nabla\mathcal E({\bf x})={\bf x}+{\bf c}-2\varphi({\bf x}). 
\end {equation}
For any given points ${\bf u},{\bf x}\in \RR^d$, consider the difference
$$
\mathcal E({\bf u})-\mathcal E({\bf x})-({\bf x}+{\bf c}-2\varphi({\bf x}))\cdot ({\bf u}-{\bf x})
$$
$$
=({\bf u}-\varphi({\bf u}))^2-({\bf x}-\varphi({\bf x}))^2+\frac 12 ({\bf x}-{\bf c})^2-\frac 12 ({\bf u}-{\bf c})^2-({\bf x}+{\bf c}-2\varphi({\bf x}))\cdot ({\bf u}-{\bf x})
$$
$$
=(({\bf u}-{\bf x})+\varphi({\bf x})-\varphi({\bf u}))\cdot ({\bf u}-\varphi({\bf u})+{\bf x}-\varphi({\bf x}))
$$
$$
-\(\frac 12 {\bf x}+\frac 12 {\bf u}-{\bf c}\)\cdot ({\bf u}-{\bf x})-({\bf x}+{\bf c}-2\varphi({\bf x}))\cdot ({\bf u}-{\bf x})
$$
$$
=(\varphi({\bf x})-\varphi({\bf u}))\cdot ({\bf u}-\varphi({\bf u})+{\bf x}-\varphi({\bf x}))+({\bf u}-{\bf x})\cdot \(\frac 12 ({\bf u}-{\bf x})+\varphi({\bf x})-\varphi({\bf u})\).
$$
Since $\varphi$ maps every point in $\RR^d$ into the closest point in $U$, by Lemma \ref {distconv} we have
$$
({\bf u}-\varphi({\bf u}))\cdot (\varphi({\bf x})-\varphi({\bf u}))\leq 0\ \ \ {\rm and}\ \ \ (\varphi({\bf x})-{\bf x})\cdot (\varphi({\bf x})-\varphi({\bf u}))\leq 0.
$$
In view of Lemma \ref {phiLip}, we have 
$$
({\bf u}-\varphi({\bf u}))\cdot (\varphi({\bf x})-\varphi({\bf u}))+(\varphi({\bf x})-{\bf x})\cdot (\varphi({\bf x})-\varphi({\bf u}))
$$
$$
=({\bf u}-{\bf x})\cdot (\varphi({\bf x})-\varphi({\bf u}))+(\varphi({\bf x})-\varphi({\bf u}))^2=O(({\bf u}-{\bf x})^2),\ \ \ {\bf u}\to {\bf x}.
$$
Thus, the expressions $({\bf u}-\varphi({\bf u}))\cdot (\varphi({\bf x})-\varphi({\bf u}))$ and $(\varphi({\bf x})-{\bf x})\cdot (\varphi({\bf x})-\varphi({\bf u}))$ have the same sign and their sum is $O(({\bf u}-{\bf x})^2)$. Then each of these expressions is $O(({\bf u}-{\bf x})^2)$, ${\bf u}\to {\bf x}$, and we finally obtain (using again Lemma \ref {phiLip})
$$
\mathcal E({\bf u})-\mathcal E({\bf x})-({\bf x}+{\bf c}-2\varphi({\bf x}))\cdot ({\bf u}-{\bf x})
$$
$$
=({\bf u}-\varphi({\bf u}))\cdot (\varphi({\bf x})-\varphi({\bf u}))-(\varphi({\bf x})-{\bf x})\cdot (\varphi({\bf x})-\varphi({\bf u}))
$$
$$
+\frac 12 ({\bf u}-{\bf x})^2+({\bf u}-{\bf x})\cdot (\varphi ({\bf x})-\varphi({\bf u}))=O(({\bf u}-{\bf x})^2), \ \ \ {\bf u}\to{\bf x}.
$$
Then by definition of the gradient, we obtain \eqref {grad_e}. By Lemma \ref {phiLip} the function $\varphi$ is continuous. Consequently, $\mathcal E\in C^1(\RR^d)$. Notice also that 
$$
\left|{\bf v}_i-{\bf u}_i\right|={\rm dist}({\bf v}_i,B[{\bf c},R/2])={\rm dist}({\bf v}_i,U),
$$
and hence,
$\varphi ({\bf v}_i)={\bf u}_i$, $i=0,1,\ldots,d$. Then 
$$
\mathcal E({\bf v}_i)=\frac {R^2}{4}+\left|{\bf v}_i-{\bf u}_i\right|^2-\frac {1}{2}\left|{\bf c}-{\bf v}_i\right|^2=\frac {1}{4}\left|{\bf v}_i-{\bf c}\right|^2-\frac {R^2}{4}=0,
$$
and $\nabla \mathcal E({\bf v}_i)={\bf v}_i+{\bf c}-2{\bf u}_i={\bf 0}$, $i=0,1,\ldots,d$.

To complete the proof of Theorem \ref {Th2} we need to show \eqref {gradient}. In view of Lemma \ref {distconv}, for any ${\bf x},{\bf y}\in \RR^d$, inequalities \eqref {x_phi} hold. Then from \eqref {grad_e} we have
$$
\left(\nabla \mathcal E({\bf x})-\nabla \mathcal E({\bf y})\right)^2=({\bf x}-{\bf y}+2(\varphi({\bf y})-\varphi({\bf x})))^2
$$
$$
=({\bf x}-{\bf y})^2+4({\bf x}-{\bf y})\cdot (\varphi({\bf y})-\varphi({\bf x}))+4(\varphi({\bf y})-\varphi({\bf x}))^2
$$
$$
=({\bf x}-{\bf y})^2+4({\bf x}-\varphi({\bf x}))\cdot (\varphi({\bf y})-\varphi({\bf x}))
+4({\bf y}-\varphi({\bf y}))\cdot (\varphi({\bf x})-\varphi({\bf y}))\leq ({\bf x}-{\bf y})^2
$$
and inequality \eqref {gradient} follows. The fact that $\mathcal E$ belongs to the class $W^2(\RR^d)$ now follows from Proposition~\ref {W2}. \hfill $\square$

\section {Proof of Theorem \ref {Th1}}\label {S9_4.1}

The following Taylor formula with the remainder in the integral form will be used in the proofs:
\begin {equation}\label {Taylor}
g(t)=g(0)+g'(0)t+\int\limits_{0}^{t}(t-u)g''(u)du,\ \ \ t\in [0,a],
\end {equation}
where $g$ is arbitrary function from $W^2_\infty[0,a]$.

We start by obtaining an upper estimate for the error of algorithm \eqref {rec1} at every point ${\bf x}\in P$. Denote 
$$
q_i({\bf x}):={\bf x}-\varphi({\bf x})+{\bf u}_i={\bf x}-\varphi({\bf x})+\frac 12 \({\bf c}+{\bf v}_i\), \ \ \ i=0,1,\ldots,d.
$$
Let $I\subset \Sigma$ be the subset such that ${\bf x}\in \varphi^{-1}(G_I)$ (the index set $I$ is defined uniquely since the sets $\varphi^{-1}(G_I)$ are pairwise disjoint). Since $\varphi({\bf x})\in G_I$, we have $b_i(\varphi({\bf x}))=0$, $i\notin I$, and $\sum_{i\in I}{b_i(\varphi({\bf x}))}=1$. For arbitrary function $f\in W^2(P)$, we will have 
$$
p_f({\bf x})-f({\bf x})=\sum\limits_{i=0}^{d}b_i(\varphi({\bf x}))\(f({\bf v}_i)+\nabla f({\bf v}_i)\cdot \(q_i({\bf x})-{\bf v}_i\)\)-f(\bf x)
$$
\begin {equation}\label {q3'}
=\sum\limits_{i\in I}b_i(\varphi({\bf x}))\(f({\bf v}_i)+\nabla f({\bf v}_i)\cdot \(q_i({\bf x})-{\bf v}_i\)-f(\bf x)\).
\end {equation}

By Definition \ref {D3.1}, for every $i\in I$, $q_i({\bf x})\in P$ and hence, the segments that join $q_i({\bf x})$ with ${\bf x}$ and with ${\bf v}_i$ are also contained in $P$. This allows us to define the functions $g_i(t):=f\({\bf x}+\frac {t}{\tau_i}(q_i({\bf x})-{\bf x})\)$, $t\in [0,\tau_i]$, where $\tau_i:=\left|q_i({\bf x})-{\bf x}\right|$, and $h_i(t):=f\({\bf v}_i+\frac {t}{\rho_i}(q_i({\bf x})-{\bf v}_i)\)$, $t\in [0,\rho_i]$, where $\rho_i=\left|q_i({\bf x})-{\bf v}_i\right|$, $i\in I$. Here we must assume that $\tau_i>0$ and $\rho_i>0$, $i\in I$. 
Then $f({\bf x})=g_i(0)$, $f({\bf v}_i)=h_i(0)$, and $\nabla f({\bf v}_i)\cdot \(q_i({\bf x})-{\bf v}_i\)=h_i'(0)\rho_i$, $i\in I$. If $\tau_i=0$, we let $g_i(0):=f({\bf x})$ and for our technique to work, we define $g'_i(0)$ and $g''_i(0)$ to be some finite numbers. If $\rho_i=0$, we set $h_i(0):=f({\bf v}_i)$ and define $h'_i(0)$ and $h_i''(0)$ to be some finite numbers. Then from \eqref {q3'} we will obtain
$$
p_f({\bf x})-f({\bf x})=\sum\limits_{i\in I}b_i(\varphi({\bf x}))\(h_i(0)+h_i'(0)\rho_i-g_i(0)\).
$$
By definition of the class $W^2(P)$, we have $g_i\in W^2_\infty[0,\tau_i]$ (if $\tau_i>0$) and $h_i\in W^2_\infty[0,\rho_i]$ (if $\rho_i>0$), $i\in I$. Since $g_i(\tau_i)=h_i(\rho_i)=f(q_i({\bf x}))$ and $g_i'(0)\tau_i=\nabla f({\bf x})\cdot (q_i({\bf x})-{\bf x})$, applying formula \eqref {Taylor} with $t=\tau_i$ and $t=\rho_i$, we will have (these equalities still hold if $\tau_i=0$ or $\rho_i=0$)
$$
p_f({\bf x})-f({\bf x})=
$$
$$
=\sum\limits_{i\in I}b_i(\varphi({\bf x}))\(h_i(\rho_i)-\int\limits_{0}^{\rho_i}{(\rho_i-u)h_i''(u)du}-g_i(\tau_i)+g_i'(0)\tau_i+\int\limits_{0}^{\tau_i}(\tau_i-u)g_i''(u)du\)
$$
$$
=\sum\limits_{i\in I}b_i(\varphi({\bf x}))\nabla f({\bf x})\cdot (q_i({\bf x})-{\bf x})
$$
$$
+\sum\limits_{i\in I}b_i(\varphi({\bf x}))\(\int\limits_{0}^{\tau_i}(\tau_i-u)g_i''(u)du-\int\limits_{0}^{\rho_i}{(\rho_i-u)h_i''(u)du}\).
$$
In view of the equalities
$$
\sum\limits_{i=0}^d{b_i(\varphi({\bf x}))q_i({\bf x})}=\sum\limits_{i=0}^{d}{b_i(\varphi({\bf x}))({\bf x}-\varphi({\bf x}))}
$$
$$
+\sum\limits_{i=0}^{d}{b_i(\varphi({\bf x})){\bf u}_i}={\bf x}-\varphi({\bf x})+\varphi({\bf x})={\bf x},
$$
we get that
$$
\sum\limits_{i\in I}b_i(\varphi({\bf x}))\nabla f({\bf x})\cdot (q_i({\bf x})-{\bf x})=\sum\limits_{i=0}^{d}{b_i(\varphi({\bf x}))\nabla f({\bf x})\cdot (q_i({\bf x})-{\bf x})}
$$
\begin {equation}\label {grad_f}
=\nabla f({\bf x})\cdot \(\sum\limits_{i=0}^{d}{b_i(\varphi({\bf x}))q_i({\bf x})}-\sum\limits_{i=0}^{d}{b_i(\varphi({\bf x})){\bf x}}\)=\nabla f({\bf x})\cdot ({\bf x}-{\bf x})=0.
\end {equation}
Then
$$
p_f({\bf x})-f({\bf x})=\sum\limits_{i\in I}b_i(\varphi({\bf x}))\(\int\limits_{0}^{\tau_i}(\tau_i-u)g_i''(u)du-\int\limits_{0}^{\rho_i}{(\rho_i-u)h_i''(u)du}\).
$$
Since $\left|g_i''(t)\right|\leq 1$ and $\left|h_i''(t)\right|\leq 1$, for almost every $t$, we obtain
$$
\left|p_f({\bf x})-f({\bf x})\right|\leq \sum\limits_{i\in I}b_i(\varphi({\bf x}))\(\int\limits_{0}^{\tau_i}\left|\tau_i-u\right| \left|g_i''(u)\right|du+\int\limits_{0}^{\rho_i}{\left|\rho_i-u\right|\left|h_i''(u)\right|du}\)
$$
$$
\leq \sum\limits_{i\in I}b_i(\varphi({\bf x}))\(\int\limits_{0}^{\tau_i}\left|\tau_i-u\right| du+\int\limits_{0}^{\rho_i}{\left|\rho_i-u\right|du}\)=\frac 12 \sum\limits_{i\in I}b_i(\varphi({\bf x}))(\tau_i^2+\rho_i^2)
$$
\begin {equation}\label {upper}
=\frac 12 \sum\limits_{i\in I}b_i(\varphi({\bf x}))(\left|q_i({\bf x})-{\bf x}\right|^2+\left|q_i({\bf x})-{\bf v}_i\right|^2).
\end {equation}
Recalling that $\left|{\bf c}-{\bf v}_i\right|=R$, $i\in I$, where $R$ is the circumradius of $T$, we will obtain
$$
\left|p_f({\bf x})-f({\bf x})\right|
$$
$$
\leq \frac 12 \sum\limits_{i\in I}b_i(\varphi({\bf x}))\(\left|\varphi({\bf x})-{\bf c}+\frac 12({\bf c}-{\bf v}_i)\right|^2+\left|{\bf x}-\varphi({\bf x})+\frac 12 ({\bf c}-{\bf v}_i)\right|^2\)
$$
$$
= \frac 12 \sum\limits_{i\in I}b_i(\varphi({\bf x}))\(\left|\varphi({\bf x})-{\bf c}\right|^2+\left|{\bf x}-\varphi({\bf x})\right|^2+({\bf c}-{\bf v}_i)\cdot ({\bf x}-{\bf c})+\frac 12 \left|{\bf c}-{\bf v}_i\right|^2\)
$$
$$
=\frac 12\left|\varphi({\bf x})-{\bf c}\right|^2+\frac 12\left|{\bf x}-\varphi({\bf x})\right|^2+\sum\limits_{i\in I}b_i(\varphi({\bf x}))({\bf c}-{\bf u}_i)\cdot ({\bf x}-{\bf c})+\frac {R^2}{4}
$$
$$
=\frac 12\left|{\bf x}-\varphi({\bf x})\right|^2+\frac 12\left|{\bf c}-\varphi({\bf x})\right|^2+({\bf c}-\varphi({\bf x}))\cdot ({\bf x}-{\bf c})+\frac {R^2}{4}
$$
$$
=\frac 12\left|{\bf x}-\varphi({\bf x})\right|^2+\frac {1}{2}\(\left|{\bf c}-\varphi({\bf x})\right|^2+2({\bf c}-\varphi({\bf x}))\cdot ({\bf x}-{\bf c})+\left|{\bf x}-{\bf c}\right|^2\)-\frac {1}{2}\left|{\bf x}-{\bf c}\right|^2+\frac {R^2}{4}
$$
$$
=\left|{\bf x}-\varphi({\bf x})\right|^2-\frac {1}{2}\left|{\bf x}-{\bf c}\right|^2+\frac {R^2}{4}=\mathcal E({\bf x}).
$$
Let $\Psi^\ast:\RR^{(d+1)^2}\to \RR$ be the mapping such that $p_f({\bf x})=\Psi^\ast(I(f))$. Then for any ${\bf x}\in P$, 
\begin {equation}\label {upp}
e_{\bf x}(W^2(P);\Psi^\ast)=\sup\limits_{f\in W^2(P)}{\left|f({\bf x})-p_f({\bf x})\right|}\leq \mathcal E({\bf x}).
\end {equation}
Since by Theorem \ref {Th2} the function $\mathcal E$ belongs to the class $W^2(\RR^d)$, we have $\mathcal E\in W^2(P)$. Hence $-\mathcal E\in W^2(P)$. By Theorem \ref {Th2} we also have $I(\mathcal E)=I(-\mathcal E)={\bf 0}$ and, in view of estimate \eqref {upp}, the function $\mathcal E$ is non-negative at every point ${\bf x}\in P$. 

Let $s(\cdot;\Psi)$ be any algorithm for recovery of $f({\bf x})$ of form \eqref {PSI}. Then for every ${\bf x}\in P$, taking into account estimate \eqref {upp} we will obtain
\begin {equation}\label {L}
\begin {split}
e_{\bf x}(W^2(P);\Psi)&=\sup\limits_{f\in W^2(P)}{\left|f({\bf x})-s(f;\Psi)\right|}\\
&\geq \max\{\left|\mathcal E({\bf x})-\Psi(I(\mathcal E))\right|,\left|-\mathcal E({\bf x})-\Psi(I(-\mathcal E))\right|\}\\
&\geq \frac 12\(\left|\mathcal E({\bf x})-\Psi({\bf 0})\right|+\left|\mathcal E({\bf x})+\Psi({\bf 0})\right|\)\\
&\geq \left|\mathcal E({\bf x})\right|=\mathcal E({\bf x})\geq e_{\bf x}(W^2(P);\Psi^\ast),
\end {split}
\end {equation}
which proves the optimality of the method $p_f({\bf x})$ for recovery of $f({\bf x})$ at every given point ${\bf x}\in P$. Since estimate \eqref {L} also holds for $\Psi=\Psi^\ast$, we obtain relation \eqref {h}.

We complete the proof of Theorem \ref {Th1} by showing the optimality of the algorithm $p_f$ for global recovery of $W^2(P)$ in any monotone seminorm $\|\ \!\cdot\ \!\|$. Let $\Phi^\ast:\RR^{(d+1)^2}\to C(P)$ be the mapping such that $p_f=\Phi^\ast(I(f))$. In view of \eqref {upp} the inequality $\left|f-p_f\right|\leq \mathcal E$ holds on $P$ for every function $f\in W^2(P)$. Then by the monotonicity of the seminorm, we have $\|f-p_f\|\leq \|\mathcal E\|$ and
\begin {equation}\label {a1}
R(W^2(P);\Phi^\ast,\|\cdot\|)=\sup\limits_{f\in W^2(P)}{\|f-\Phi^\ast(I(f))\|}\leq \|\mathcal E\|.
\end {equation}
Let now $S(\cdot;\Phi)$ be any global recovery algorithm of form \eqref {alg2}. Then, by an argument similar to \eqref {L}, we obtain
$$
R(W^2(P);\Phi,\|\cdot\|)=\sup\limits_{f\in W^2(P)}{\|f-\Phi(I(f))\|}
$$
$$
\geq \max\{\|\mathcal E-\Phi(I(\mathcal E))\|,\|-\mathcal E-\Phi(I(-\mathcal E))\|\}
$$
$$
\geq \frac 12\(\|\mathcal E-\Phi({\bf 0})\|+\|\mathcal E+\Phi({\bf 0})\|\)\geq \|\mathcal E\|\geq R(W^2(P);\Phi^\ast,\|\cdot\|),
$$
which shows the optimality of the algorithm $p_f$ on $W^2(P)$ and if we let $\Phi=\Phi^\ast$ in the above estimate, we will obtain  \eqref {h1}.\hfill $\square$

\medskip

{\bf Remark.} Corollary \ref {q1'} can also be derived from the proof of Theorem \ref {Th1}. Indeed, since $p_f({\bf x})=0$ for every function $f\in W^2(P)$ with $I(f)={\bf 0}$, relation \eqref {upp} implies that
\begin {equation}\label {j1'}
\sup\limits_{f\in W^2(P)\atop I(f)={\bf 0}}{f({\bf x})}\leq \mathcal E({\bf x}).
\end {equation}
Since by Theorem \ref {Th2}, the function $\mathcal E$ belongs to the class $W^2(P)$ and $I(\mathcal E)={\bf 0}$, we get equality in \eqref {j1'}.

\section {Analysis of the function $\mathcal E$}\label {S10}

{\bf Proof of Proposition \ref {partition}.} We first show that the sets $G_I$ and $F_I$ are mutually orthogonal. Let $i,j\in I$ and ${\bf y},{\bf z}\in F_I$. Then the equalities $\left|{\bf y}-{\bf v}_i\right|^2=\left|{\bf y}-{\bf v}_j\right|^2$ and $\left|{\bf z}-{\bf v}_i\right|^2=\left|{\bf z}-{\bf v}_j\right|^2$ imply that $2{\bf y}\cdot ({\bf v}_j-{\bf v}_i)={\bf v}_j^2-{\bf v}_i^2=2{\bf z}\cdot ({\bf v}_j-{\bf v}_i)$ and hence, $({\bf y}-{\bf z})\cdot ({\bf v}_i-{\bf v}_j)=0$, $i,j\in I$. Consequently, $({\bf y}-{\bf z})\cdot ({\bf u}_i-{\bf u}_j)=0$, $i,j\in I$. Thus, $G_I$ is mutually orthogonal with $F_I$ and hence, with $(1/2)(F_I-{\bf c})$.

We next show \eqref {p2}. Let ${\bf x}\in \varphi^{-1}(G_I)$. Show that $2({\bf x}-\varphi({\bf x}))+{\bf c}\in F_I$. Since $\varphi({\bf x})$ is the point in $U$ closest to ${\bf x}$, by Lemma \ref {distconv} we have 
\begin {equation}\label {E1'}
({\bf x}-\varphi({\bf x}))\cdot ({\bf u}_j-\varphi({\bf x}))\leq 0, \ \ \ j\in \Sigma. 
\end {equation}
Since $b_i(\varphi({\bf x}))>0$, $i\in I$, and
$$
\sum\limits_{i\in I}{b_i(\varphi({\bf x}))({\bf x}-\varphi({\bf x}))\cdot ({\bf u}_i-\varphi({\bf x}))}=({\bf x}-\varphi({\bf x}))\cdot \(\sum\limits_{i\in I}{b_i(\varphi({\bf x}))}({\bf u}_i-\varphi({\bf x}))\)
$$ 
$$
=({\bf x}-\varphi({\bf x}))\cdot (\varphi({\bf x})-\varphi({\bf x}))=0,
$$
we have $({\bf x}-\varphi({\bf x}))\cdot ({\bf u}_i-\varphi({\bf x}))= 0$, $i\in I$.
Since $\left|{\bf c}-{\bf v}_i\right|=\left|{\bf c}-{\bf v}_j\right|$, $i,j\in \Sigma$, for every $i\in I$ and $j\in \Sigma$, we have
$$
\left|2({\bf x}-\varphi({\bf x}))+{\bf c}-{\bf v}_i\right|^2-\left|2({\bf x}-\varphi({\bf x}))+{\bf c}-{\bf v}_j\right|^2
$$
$$
=4({\bf x}-\varphi({\bf x}))\cdot ({\bf v}_j-{\bf v}_i)=8({\bf x}-\varphi({\bf x}))\cdot ({\bf u}_j-{\bf u}_i)
$$
$$
=8({\bf x}-\varphi({\bf x}))\cdot ({\bf u}_j-\varphi({\bf x}))+8({\bf x}-\varphi({\bf x}))\cdot (\varphi({\bf x})-{\bf u}_i)=8({\bf x}-\varphi({\bf x}))\cdot ({\bf u}_j-\varphi({\bf x}))\leq 0.
$$
Consequently, $2({\bf x}-\varphi({\bf x}))+{\bf c}\in W_i$ for every $i\in I$. Then we have $2({\bf x}-\varphi({\bf x}))+{\bf c}\in F_I$, which implies that ${\bf x}-\varphi({\bf x})\in (1/2)(F_I-{\bf c})$. Since $\varphi({\bf x})\in G_I$, we have that ${\bf x}\in G_I\oplus (1/2)(F_I-{\bf c})$.

Let now ${\bf x}\in G_I\oplus (1/2)(F_I-{\bf c})$. Then ${\bf x}={\bf u}+{\bf z}$, where ${\bf u}\in G_I$ and ${\bf z}\in  (1/2)(F_I-{\bf c})$. We will show that $\varphi({\bf x})={\bf u}\in G_I$. For every $i\in \Sigma$, we have
\begin {equation}\label {E3}
({\bf x}-{\bf u})\cdot ({\bf u}_i-{\bf u})={\bf z}\cdot ({\bf u}_i-{\bf u})=(1/2) ({\bf y}-{\bf c})\cdot ({\bf u}_i-{\bf u}),
\end {equation}
where ${\bf y}$ is some point from $F_I$. For every $i\in I$ and $j\in \Sigma$, we then have $\left|{\bf y}-{\bf v}_i\right|\leq \left|{\bf y}-{\bf v}_j\right|$ and hence,
\begin {equation}\label {25a}
\begin {split}
0\geq \left|{\bf y}-{\bf v}_i\right|^2&-\left|{\bf y}-{\bf v}_j\right|^2=\left|({\bf y}-{\bf c})+({\bf c}-{\bf v}_i)\right|^2-\left|({\bf y}-{\bf c})+({\bf c}-{\bf v}_j)\right|^2\\
&=2({\bf y}-{\bf c})\cdot ({\bf v}_j-{\bf v}_i)=4({\bf y}-{\bf c})\cdot ({\bf u}_j-{\bf u}_i)
\end {split}
\end {equation}
with equality being true whenever $j\in I$. Then for every $i,j\in I$, 
$$
({\bf y}-{\bf c})\cdot ({\bf u}_i-{\bf u})-({\bf y}-{\bf c})\cdot ({\bf u}_j-{\bf u})=({\bf y}-{\bf c})\cdot ({\bf u}_i-{\bf u}_j)=0.
$$
Consequently, the quantity $({\bf y}-{\bf c})\cdot ({\bf u}_i-{\bf u})$ has the same value for every $i\in I$, which we denote by $a$. Since ${\bf u}\in G_I$, we have $b_i({\bf u})=0$, $i\notin I$, and $\sum_{i\in I}b_i({\bf u})=1$. Then
$$
a=\sum_{i\in I}b_i({\bf u})a=\sum_{i\in I}b_i({\bf u})({\bf y}-{\bf c})\cdot ({\bf u}_i-{\bf u})=({\bf y}-{\bf c})\cdot \(\sum_{i\in I}b_i({\bf u})({\bf u}_i-{\bf u})\)=({\bf y}-{\bf c})\cdot ({\bf u}-{\bf u})=0.
$$
Hence $({\bf y}-{\bf c})\cdot ({\bf u}_i-{\bf u})=0$, $i\in I$. Furthermore, for every $j\notin I$, using any $i\in I$, we can write (see \eqref {25a})
$$
({\bf y}-{\bf c})\cdot ({\bf u}_j-{\bf u})=({\bf y}-{\bf c})\cdot ({\bf u}_i-{\bf u})+({\bf y}-{\bf c})\cdot ({\bf u}_j-{\bf u}_i)=({\bf y}-{\bf c})\cdot ({\bf u}_j-{\bf u}_i)\leq 0.
$$
Using \eqref {E3} we now obtain that $({\bf x}-{\bf u})\cdot ({\bf u}_i-{\bf u})\leq 0$, $i\in \Sigma$. If now ${\bf t}=\sum_{i=0}^{d}b_i({\bf t}){\bf u}_i$ is any point in $U$, then
$$
({\bf x}-{\bf u})\cdot ({\bf t}-{\bf u})=({\bf x}-{\bf u})\cdot \(\sum_{i=0}^{d}b_i({\bf t}){\bf u}_i-{\bf u}\)=\sum_{i=0}^{d}b_i({\bf t})({\bf x}-{\bf u})\cdot ({\bf u}_i-{\bf u})\leq 0,
$$
which in view of Lemma \ref {distconv} implies that $\varphi({\bf x})={\bf u}\in G_I$ and hence, ${\bf x}\in \varphi^{-1}(G_I)$. This completes the proof of the first equality in \eqref {p2}. The sets $D_I$ and $F_I$ are mutually orthogonal because $G_I$ and $F_I$ are. Then the second equality in \eqref {p2} follows immediately.

It remains to show that $\varphi^{-1}(G_I)$ is convex and has a non-empty interior. It is not difficult to see that the set $G_I$ is convex and $F_I$ is convex as the intersection of convex sets $W_i$, $i\in I$. Then $(1/2)(F_I-{\bf c})$ is also convex. By \eqref {p2}, if ${\bf x}_1,{\bf x}_2\in \varphi^{-1}(G_I)$, then ${\bf x}_i={\bf y}_i+{\bf z}_i$, where ${\bf y}_i\in G_I$ and ${\bf z}_i\in (1/2)(F_I-{\bf c})$, $i=1,2$, and for every $t\in [0,1]$, we have $t{\bf x}_1+(1-t){\bf x}_2=t{\bf y}_1+(1-t){\bf y}_2+t{\bf z}_1+(1-t){\bf z}_2\in G_I\oplus (1/2)(F_I-{\bf c})=\varphi^{-1}(G_I)$, which shows the convexity of $\varphi^{-1}(G_I)$. 

We next show that the set
$$
F_I^\circ :=\bigcap _{i\in I} W_i\setminus \bigcup_{i\in \Sigma\setminus I} W_i
$$
is non-empty. To do this we first show that $F_{\Sigma\setminus \{i\}}^\circ$ is non-empty for every $i\in \Sigma$. Let ${\bf n}\cdot {\bf x}=a$ be the equation of the hyperplane $P_i$ that contains the face $D_{\Sigma\setminus \{i\}}$ of $T$ (${\bf n}$ is chosen to be a unit vector). Let ${\bf c}_i$ be the circumcenter of $D_{\Sigma\setminus \{i\}}$ and let $R_i$ be its circumradius. Let ${\bf n}\cdot {\bf v}_i=b$ and ${\bf y}:={\bf c}_i+t{\bf n}$ (notice that $b\neq a$). For every $j\in \Sigma\setminus \{i\}$, since ${\bf n}\cdot {\bf c}_i={\bf n}\cdot {\bf v}_j=a$, we have ${\bf c}_i-{\bf v}_j\ \bot \ {\bf n}$. Then $({\bf y}-{\bf v}_j)^2=({\bf c}_i-{\bf v}_j)^2+t^2=R_i^2+t^2$. Consequently, 
$$
({\bf y}-{\bf v}_i)^2=({\bf c}_i-{\bf v}_i)^2+2t({\bf c}_i-{\bf v}_i)\cdot {\bf n}+t^2=({\bf c}_i-{\bf v}_i)^2+2t(a-b)+t^2,
$$
which is strictly greater than $R_i^2+t^2$ if and only if $2t(a-b)>R_i^2-({\bf c}_i-{\bf v}_i)^2$. Then there is $t$ such that $\left|{\bf y}-{\bf v}_j\right|=\sqrt{R_i^2+t^2}$, $j\in \Sigma\setminus \{i\}$, and $\sqrt {R_i^2+t^2}<\left|{\bf y}-{\bf v}_i\right|$. Then ${\bf y}\in W_j$, $j\in \Sigma\setminus \{i\}$, and ${\bf y}\notin W_i$. This shows that ${\bf y}\in F_{\Sigma\setminus \{i\}}^\circ$ and, hence, the set $F_{\Sigma\setminus \{i\}}^\circ$ is non-empty.

Let now $I\subset \Sigma$ be any non-empty proper subset. Show that $F_I^\circ\neq \emptyset$. Let ${\bf z}:=\sum\limits_{i\in \Sigma\setminus I}\alpha_i {\bf y}_i$, where ${\bf y}_i\in F_{\Sigma\setminus \{i\}}^\circ$ and $\alpha_i$'s are strictly positive numbers that add up to $1$.  

Fix $j\in I$. Recall that $X=\{{\bf v}_0,\ldots,{\bf v}_d\}$. Then for every $k\in \Sigma\setminus I$, we will have ${\bf y}_k\notin W_k$ and ${\bf y}_k\in W_j$. Then $\left|{\bf y}_k-{\bf v}_k\right|>{\rm dist}({\bf y}_k,X)=\left|{\bf y}_k-{\bf v}_j\right|$. For every $i\in \Sigma\setminus I$, $i\neq k$, we have ${\bf y}_i\in W_k$ and ${\bf y}_i\in W_j$. Then $\left|{\bf y}_i-{\bf v}_k\right|={\rm dist}({\bf y}_i,X)=\left|{\bf y}_i-{\bf v}_j\right|$. Consequently, 
$$
({\bf z}-{\bf v}_k)^2-({\bf z}-{\bf v}_j)^2=2{\bf z}\cdot ({\bf v}_j-{\bf v}_k)+{\bf v}_k^2-{\bf v}_j^2
$$
$$
=\sum\limits_{i\in \Sigma\setminus I}\alpha_i(2{\bf y}_i\cdot ({\bf v}_j-{\bf v}_k)+{\bf v}_k^2-{\bf v}_j^2)=\sum\limits_{i\in \Sigma\setminus I}\alpha_i\(({\bf y}_i-{\bf v}_k)^2-({\bf y}_i-{\bf v}_j)^2\)
$$
$$
=\alpha_k\(({\bf y}_k-{\bf v}_k)^2-({\bf y}_k-{\bf v}_j)^2\)>0.
$$
Then ${\bf z}\notin W_k$, $k\in \Sigma\setminus I$.
When $k\in I$ and $i\in \Sigma\setminus I$, we will have $\left|{\bf y}_i-{\bf v}_k\right|={\rm dist}({\bf y}_i,X)=\left|{\bf y}_i-{\bf v}_j\right|$. Running the same argument we will obtain $({\bf z}-{\bf v}_k)^2-({\bf z}-{\bf v}_j)^2=0$. Then ${\bf z}\in W_j$, $j\in I$, and consequently, ${\bf z}\in F_I^\circ$. Thus, $F_I^\circ\neq \emptyset$.

If $I=\Sigma$, then $\varphi^{-1}(G_I)$ equals the interior of the simplex $U$, and thus, the interior of $\varphi^{-1}(G_\Sigma)$ is non-empty. If $I$ is a non-empty proper subset of $\Sigma$, then the set $D_I\oplus F^\circ_I$ is non-empty because $D_I$ and $F^\circ_I$ are. By Lemma \ref {L10}, the direct sum $D_I\oplus F_I^\circ$ is open in $\RR^d$ and hence, the interiors of $D_I\oplus F_I$ and of $\varphi^{-1}(G_I)=\frac 12 (D_I\oplus F_I)$ are non-empty.\hfill $\square$

\medskip

For a non-empty subset $I\subset \Sigma$, we define
\begin {equation}\label {WQ}
\W Q_I:=\left\{\sum_{i\in I}c_i{\bf u}_i : \sum_{i\in I}c_i=1\right\}=\frac {1}{2}\({\bf c}+Q_I\). 
\end {equation}

\medskip

{\bf Proof of Theorem \ref {E}.} Let $\W {\bf c}_I:=\frac 12 ({\bf c}+{\bf c}_I)$. Then $\W {\bf c}_I$ is the circumcenter of the face $G_I$ of the simplex $U$. Notice that ${\bf c},{\bf c}_I\in H_I$. Then by Lemma \ref {L10a}, we have 
\begin {equation}\label {qq}
({\bf c}-{\bf c}_I)\cdot ({\bf a}-{\bf b})=0 \ \ \ {\rm for \ any}\ \ \  {\bf a},{\bf b}\in Q_I. 
\end {equation}
The affine subspace $H_I$ is also mutually orthogonal with the affine subspace $\W Q_I$.

Let ${\bf x}\in \varphi ^{-1}(G_I)$ be arbitrary. It is not difficult to verify that $\varphi({\bf x})$ is the orthogonal projection of ${\bf x}$ onto $\W Q_I$. Consequently,
\begin {equation}\label {5w}
({\bf x}-\varphi({\bf x}))\cdot ({\bf p}-{\bf q})=0,\ \ \ {\bf p},{\bf q}\in \W Q_I.
\end {equation}
Let ${\bf z}:=\W {\bf c}_I+{\bf x}-\varphi({\bf x})$. Since $\varphi({\bf x}),\W {\bf c}_I\in \W Q_I$, for every ${\bf v}\in H_I$, we have
$$
({\bf x}-{\bf z})\cdot ({\bf v}-{\bf z})=\(\varphi({\bf x})-\W {\bf c}_I\)\cdot ({\bf v}-\W {\bf c}_I-({\bf x}-\varphi({\bf x})))=(\varphi({\bf x})-\W {\bf c}_I)\cdot ({\bf v}-\W {\bf c}_I).
$$
Since $\W {\bf c}_I\in H_I$ and $H_I$ is mutually orthogonal with $\W Q_I$, we have $(\varphi({\bf x})-\W {\bf c}_I)\cdot ({\bf v}-\W {\bf c}_I)=0$. Then $({\bf x}-{\bf z})\cdot ({\bf v}-{\bf z})=0$. Observe also that ${\bf z}\in H_I$. Indeed, for every $i,j\in I$, using \eqref {qq} and \eqref {5w}, we have
$$
({\bf z}-{\bf c}_I)\cdot ({\bf v }_i-{\bf v}_j)=\frac {1}{2}\({\bf c}-{\bf c}_I\)\cdot ({\bf v }_i-{\bf v}_j)+2({\bf x}-\varphi({\bf x}))\cdot ({\bf u}_i-{\bf u}_j)=0.
$$
Now an argument similar to \eqref {5E} (see the Appendix) shows that $\left|{\bf z}-{\bf v}_i\right|=\left|{\bf z}-{\bf v}_j\right|$, $i,j\in I$. Then ${\bf z}\in H_I$. In Lemma \ref {distconv} it is sufficient to assume that $K$ is only closed and convex. Then applying Lemma \ref {distconv} with $K=H_I$ we obtain that ${\bf z}$ is the point in $H_I$ closest to ${\bf x}$.

Denote ${\bf q}:=\varphi({\bf x})-{\bf c}+\W {\bf c}_I=\varphi({\bf x})+\frac {1}{2}\({\bf c}_I-{\bf c}\)$. Since $\varphi({\bf x})\in \W Q_I$, we have ${\bf q}\in Q_I$. Then for every ${\bf w}\in Q_I$, using \eqref {qq} and the fact that $Q_I=2\W Q_I -{\bf c}_I$, we will obtain
$$
({\bf x}-{\bf q})\cdot ({\bf w}-{\bf q})=({\bf x}-\varphi({\bf x}))\cdot ({\bf w}-{\bf q})-\frac {1}{2}({\bf c}_I-{\bf c})\cdot ({\bf w}-{\bf q})=2({\bf x}-\varphi({\bf x}))\cdot ({\bf w}_1-{\bf q}_1)
$$
for some ${\bf w}_1,{\bf q}_1\in \W Q_I$. Then  \eqref {5w} implies that $({\bf x}-{\bf q})\cdot ({\bf w}-{\bf q})=0$.
By Lemma \ref {distconv} (which still holds if the compactness assumption on $K$ is replaced by the assumption that $K$ is closed), ${\bf q}$ is the point in $Q_I$ closest to ${\bf x}$. Consequently,
$$
{\rm dist}({\bf x},H_I)=\left|{\bf x}-{\bf z}\right|=\left|\varphi({\bf x})-\W {\bf c}_I\right|\ \ \ {\rm and}\ \ \ {\rm dist}({\bf x},Q_I)=\left|{\bf x}-\varphi({\bf x})+{\bf c}-\W {\bf c}_I\right|,\ \ \ {\bf x}\in \varphi^{-1}(G_I).
$$ 
In view of \eqref {qq}, we have $R^2=R_I^2+\left|{\bf c}-{\bf c}_I\right|^2=R_I^2+4\left|{\bf c}-\W {\bf c}_I\right|^2$ and in view of \eqref {5w}, \eqref {WQ}, and \eqref {qq}, we have $({\bf x}-\varphi({\bf x}))\cdot(\varphi({\bf x})-\W {\bf c}_I)=({\bf c}-\W {\bf c}_I)\cdot (\varphi({\bf x})-\W {\bf c}_I)=0$. Then for every ${\bf x}\in \varphi^{-1}(G_I)$, we obtain
$$
\mathcal E({\bf x})=\frac {R^2}{4}+\left|{\bf x}-\varphi({\bf x})\right|^2-\frac {1}{2}\left|{\bf x}-{\bf c}\right|^2
$$
$$
=\frac {R_I^2}{4}+\left|{\bf c}-\W {\bf c}_I\right|^2+\left|{\bf x}-\varphi({\bf x})\right|^2-\frac {1}{2}\left|({\bf x}-\varphi({\bf x}))+(\varphi({\bf x})-\W {\bf c}_I)-({\bf c}-\W {\bf c}_I)\right|^2
$$
$$
=\frac {R_I^2}{4}+\frac 12 \left|{\bf c}-\W {\bf c}_I\right|^2+\frac 12 \left|{\bf x}-\varphi({\bf x})\right|^2-\frac 12 \left|\varphi({\bf x})-\W {\bf c}_I\right|^2+({\bf x}-\varphi({\bf x}))\cdot ({\bf c}-\W {\bf c}_I)
$$
$$
=\frac {R_I^2}{4}-\frac 12 \left|\varphi({\bf x})-\W {\bf c}_I\right|^2+\frac 12 \left|{\bf x}-\varphi({\bf x})+{\bf c}-\W {\bf c}_I\right|^2
$$
$$
=\frac {R_I^2}{4}-\frac 12 \({\rm dist}({\bf x},H_I)\)^2+\frac 12 \({\rm dist}({\bf x},Q_I)\)^2.
$$ 
By continuity, equality \eqref {E1} holds on the closure $\OL {\varphi^{-1}(G_I)}$.

To show equality \eqref {E2}, we recall that $\nabla \mathcal E({\bf x})={\bf x}+{\bf c}-2\varphi(\bf x)$, ${\bf x}\in \RR^d$, see Theorem \ref {Th2}. Let ${\bf x},{\bf y}\in \varphi^{-1}(G_I)$ be arbitrary. Since $\varphi({\bf x})$ and $\varphi({\bf y})$ are both contained in $G_I\subset W Q_I$, the vectors ${\bf x}-\varphi({\bf x})$ and ${\bf y}-\varphi({\bf y})$ are both orthogonal to the vector $\varphi({\bf x})-\varphi({\bf y})$, see \eqref {5w}. Then
$$
\left(\nabla \mathcal E({\bf x})-\nabla \mathcal E({\bf y})\right)^2=({\bf x}-{\bf y}-2(\varphi({\bf x})-\varphi({\bf y})))^2
$$
$$
=({\bf x}-{\bf y})^2-4({\bf x}-\varphi({\bf x}))\cdot (\varphi({\bf x})-\varphi({\bf y}))
+4({\bf y}-\varphi({\bf y}))\cdot (\varphi({\bf x})-\varphi({\bf y}))= ({\bf x}-{\bf y})^2.
$$
By continuity, equality \eqref {E2} extends to the closure of $\varphi^{-1}(G_I)$.
 \hfill $\square$

\medskip

To prove Theorem \ref {nd} we will need the following lemma.

\begin {lemma}\label {phi_x}
Let $T$ be a non-degenerate simplex in $\RR^d$ with vertices ${\bf v}_0,\ldots,{\bf v}_d$. Let also $i\in \Sigma$, $I:=\Sigma\setminus \{i\}$, and assume that the circumcenter ${\bf c}$ and the vertex ${\bf v}_i$ of $T$ lie in the same closed half-space relative to the hyperplane $Q_I$. Let ${\bf n}$ be a unit normal vector to the facet $T_I$ of $T$. Then for every point ${\bf x}$ in the relative interior $D_I$ of $T_I$, we have $\frac {\partial \mathcal E}{\partial {\bf n}}({\bf x})=0$ if and only if $\varphi({\bf x})\in \OL G_I$.
\end {lemma}
\begin {proof}
For every ${\bf x}\in D_I$, we have $\W {\bf x}:=\frac {1}{2}({\bf x}+{\bf c})\in G_I$. The hyperplane $\W Q_I$ that contains the facet $\OL G_I$ of the simplex $U$ has equation $({\bf t}-\W {\bf x})\cdot {\bf n}=0$. By Theorem \ref {Th2}, we have
$$
\frac {\partial \mathcal E}{\partial {\bf n}}({\bf x})=\nabla \mathcal E({\bf x})\cdot {\bf n}=({\bf x}+{\bf c}-2\varphi({\bf x}))\cdot {\bf n}=2(\W {\bf x}-\varphi({\bf x}))\cdot {\bf n}.
$$
Then $\frac {\partial \mathcal E}{\partial {\bf n}}({\bf x})=0$ if and only if $\varphi({\bf x})\in \W Q_I$. Since $\varphi({\bf x})\in U$, this is equivalent to $\varphi({\bf x})\in \OL G_I$.
\end {proof}

{\bf Proof of Theorem \ref {nd}.} Let $i\in \Sigma$ be arbitrary, $I:=\Sigma\setminus \{i\}$, and let ${\bf x}$ be an arbitrary point in $D_I$. In view of Lemma \ref {phi_x}, we need to show that $\varphi({\bf x})\in \OL G_I$. Assume to the contrary that $\varphi({\bf x})\notin \OL G_I$. Then $b_i(\varphi({\bf x}))>0$. The assumptions of the theorem imply that the orthogonal projection $\W {\bf u}_i$ of the vertex ${\bf u}_i$ of $U$ onto the hyperplane $\W Q_I$ belongs to~$\OL G_I$. 

Let ${\bf z}$ be the orthogonal projection of the point $\varphi({\bf x})$ onto the hyperplane $\W Q_I$.
If $\varphi({\bf x})={\bf u}_i$, then ${\bf z}=\W {\bf u}_i\in \OL G_I$.
If $\varphi({\bf x})\neq {\bf u}_i$, we have
$$
{\bf y}_t:=t\varphi({\bf x})+(1-t){\bf u}_i=\sum\limits_{j=0\atop j\neq i}^{d}tb_j(\varphi({\bf x})){\bf u}_j+(tb_i(\varphi({\bf x}))+1-t){\bf u}_i.
$$
Let $t_0:=\frac {1}{1-b_i(\varphi({\bf x}))}$. Then $t_0b_i(\varphi({\bf x}))+1-t_0=0$ and we have ${\bf y}_{t_0}\in \OL G_I$. Since $\varphi({\bf x})=\frac {1}{t_0}{\bf y}_{t_0}+\(1-\frac {1}{t_0}\){\bf u}_i$, where $1/t_0\in (0,1)$, we have ${\bf z}=\frac {1}{t_0}{\bf y}_{t_0}+\(1-\frac {1}{t_0}\)\W {\bf u}_i\in \OL G_I\subset U$. The line through the points $\varphi({\bf x})$ and ${\bf z}$ is also orthogonal to the hyperplane $Q_I$. By assumption, ${\bf c}\in T$. Then ${\bf c}\in U$ and hence $b_i(c)\geq 0$. Consequently, $\varphi({\bf x})$ and ${\bf c}$ lie in the same closed half-space relative to the hyperplane $\W Q_I$. In view of \eqref {WQ}, the point ${\bf c}$ and the hyperplane $Q_I$ lie in different half-spaces relative to the hyperplane $\W Q_I$ (or ${\bf c}\in Q_I=\W Q_I$). Then $\varphi({\bf x})$ and $Q_I$ lie in different half-spaces relative to the hyperplane $\W Q_I$ (or $Q_I=\W Q_I$). Consequently, ${\bf z}\in \W Q_I$ is strictly closer to $Q_I$ than $\varphi({\bf x})$.
By the Pythagorean Theorem, the point ${\bf z}\in U$ will be strictly closer to the point ${\bf x}$ (which lies in $Q_I$) than the point $\varphi({\bf x})$. This contradicts the definition of the point $\varphi({\bf x})$. Thus, $\varphi({\bf x})\in \OL G_I$ and by Lemma \ref {phi_x}, we have $\frac {\partial \mathcal E}{\partial {\bf n}}({\bf x})=0$.

\section {Appendix}\label {A11}

{\bf Proof of Lemma \ref {distconv}.}
Let ${\bf y}$ be a point in $K$ closest to ${\bf x}$ (there exists at least one such point since $K$ is compact). We will first prove \eqref {q}. Assume to the contrary that there is a point ${\bf z}\in K$ such that $({\bf x}-{\bf y})\cdot ({\bf z}-{\bf y})>0$. Define the function $g(t):=({\bf x}-(t{\bf z}+(1-t){\bf y}))^2$, $t\in [0,1]$. Then $g'(t)=2t({\bf x}-{\bf z})^2+(2-4t)({\bf x}-{\bf z})\cdot ({\bf x}-{\bf y})+2(t-1)({\bf x}-{\bf y})^2$ and, in particular, $g'(0)=2({\bf x}-{\bf y})\cdot ({\bf y}-{\bf z})<0$. Hence, there is $t\in (0,1]$ sufficiently close to zero such that
$$
g(t)=({\bf x}-(t{\bf z}+(1-t){\bf y}))^2<g(0)=({\bf x}-{\bf y})^2.
$$
Since $t{\bf z}+(1-t){\bf y}\in K$, we get a contradiction with the fact that ${\bf y}$ is closest in $K$ to ${\bf x}$.

Let now ${\bf y}\in K$ be a point such that \eqref {q} holds. Then for every point ${\bf z}\in K$,
$$
({\bf x}-{\bf z})^2=(({\bf x}-{\bf y})-({\bf z}-{\bf y}))^2
$$
$$
=({\bf x}-{\bf y})^2-2({\bf x}-{\bf y})\cdot({\bf z}-{\bf y})+({\bf z}-{\bf y})^2\geq ({\bf x}-{\bf y})^2,
$$
which implies that ${\bf y}$ is the point in $K$ closest to ${\bf x}$.

To show the uniqueness, assume to the contrary that there is another point ${\bf u}\in K$ closest to ${\bf x}$. In view of \eqref {q}, we have
$
({\bf x}-{\bf y})\cdot ({\bf u}-{\bf y})\leq 0
$
and $
({\bf x}-{\bf u})\cdot ({\bf y}-{\bf u})\leq 0,
$
and hence,
$$
({\bf u}-{\bf y})^2=({\bf x}-{\bf y})\cdot ({\bf u}-{\bf y})+({\bf u}-{\bf x})\cdot ({\bf u}-{\bf y})\leq 0,
$$
which implies that ${\bf y}={\bf u}$. \hfill $\square$

\medskip

For the proof of Proposition \ref {partition}, we need the following two auxiliary statements.
\begin {lemma}\label {L10a}
Let $\{{\bf v}_i\}_{i\in \Sigma}$ be an affinely independent set of points in $\RR^d$. 
Let $I\subset \Sigma$ be any non-empty subset. Then the sets $Q_I$ and $H_I$ are mutually orthogonal affine subspaces of $\RR^d$ with ${\rm dim}\ \!Q_I=\#I-1$ and ${\rm dim}\ \! H_I=d+1-\#I$. Furthermore, $Q_I\cap H_I=\{c_I\}$ and $Q_I\oplus H_I=\RR^d$.
\end {lemma}
{\bf Proof.}
Fix an index $i_0\in I$. If $I\setminus \{i_0\}=\emptyset$, we have $Q_I=\{{\bf v}_{i_0}\}$, $H_I=\RR^d$, ${\bf c}_I={\bf v}_{i_0}$, and the assertion of the lemma follows trivially. Therefore, assume that $I\setminus \{i_0\}\neq \emptyset$ and let $S_I:={\rm span}\{{\bf v}_i-{\bf v}_{i_0}\}_{i\in I\setminus \{i_0\}}$. Then $Q_I={\bf v}_{i_0}+S_I$. For every point ${\bf x}\in \RR^d$ and index $i\in I\setminus \{i_0\}$, the following equalities are equivalent:
\begin {equation}\label {5E}
\begin {split}
\left|{\bf x}-{\bf v}_{i_0}\right|^2&=\left|{\bf x}-{\bf v}_{i}\right|^2,\\
\(\({\bf x}-{\bf c}_I\)+\({\bf c}_I-{\bf v}_{i_0}\)\)^2&=\(\({\bf x}-{\bf c}_I\)+\({\bf c}_I-{\bf v}_i\)\)^2,\\
\({\bf x}-{\bf c}_I\)^2+2\({\bf x}-{\bf c}_I\)\cdot \({\bf c}_I-{\bf v}_{i_0}\)+\left|{\bf c}_I-{\bf v}_{i_0}\right|^2&=\({\bf x}-{\bf c}_I\)^2+2\({\bf x}-{\bf c}_I\)\cdot \({\bf c}_I-{\bf v}_{i}\)+\left|{\bf c}_I-{\bf v}_{i}\right|^2,\\
\({\bf x}-{\bf c}_I\)\cdot \({\bf c}_I-{\bf v}_{i_0}\)&=\({\bf x}-{\bf c}_I\)\cdot \({\bf c}_I-{\bf v}_{i}\)\\
\({\bf x}-{\bf c}_I\)\cdot \({\bf v}_i-{\bf v}_{i_0}\)&=0.
\end {split}
\end {equation}
Then ${\bf x}\in H_I$ if and only if ${\bf x}-{\bf c}_I\bot \{{\bf v}_i-{\bf v}_{i_0}\}_{i\in I\setminus \{i_0\}}$. Consequently, $H_I={\bf c}_I+S_I^\bot$. This implies that the sets $Q_I$ and $H_I$ are mutually orthogonal affine subspaces. The set $Q_I$ can be rewritten as $Q_I={\bf c}_I+S_I$. Then $Q_I\cap H_I=\{{\bf c}_I\}$.
Since $\{{\bf v}_i\}_{i\in I}$ is affinely independent, the system of vectors $\{{\bf v}_i-{\bf v}_{i_0}\}_{i\in I\setminus \{i_0\}}$ is linearly independent. Then ${\rm dim}\ \!Q_I={\rm dim}\ \! S_I=\# I -1$. Since $\RR^d=S_I\oplus S_I^\bot$, we have  ${\rm dim}\ \!Q_I+{\rm dim}\ \! H_I=d$ and hence, ${\rm dim}\ \! H_I=d+1-\#I$. Furthermore, $\RR^d=Q_I\oplus H_I$.\hfill $\square$

\begin {lemma}\label {L10}
Let $\{{\bf v}_i\}_{i\in \Sigma}$ be an affinely independent set of points in $\RR^d$. 
Let $I\subset \Sigma$ be a non-empty proper subset. Then the set $D_I\oplus F_I^\circ$ is open in $\RR^d$. 
\end {lemma}
{\bf Proof.} If $I$ consists of just one index, say $i_0$, we have $D_I=\{{\bf v}_{i_0}\}$, $F_I^\circ={\rm int}\ \! W_{i_0}$, and $D_I\oplus F_I^\circ={\bf v}_{i_0}+{\rm int}\ \! W_{i_0}$, which is open in $\RR^d$. Therefore, assume that $\# I\geq 2$. Observe that 
$$
F^\circ _I=\{{\bf x}\in H_I : \left|{\bf x}-{\bf v}_k\right|>\left|{\bf x}-{\bf v}_j\right|,\ k\in \Sigma\setminus I,\ j\in I\},
$$
which is clearly open relative to $H_I$. For a point ${\bf x}\in \RR^d$, let $\{\beta_i({\bf x})\}_{i=0}^d$ be the unique collection of numbers such that $\sum_{i=0}^d\beta_i({\bf x})=1$ and ${\bf x}=\sum_{i=0}^d\beta_i({\bf x}){\bf v}_i$. Then $Q_I=\{{\bf x}\in \RR^d : \beta_i({\bf x})=0,\ i\in \Sigma\setminus I\}$ and
$$
D_I=\{{\bf x}\in Q_I : \beta_i({\bf x})>0,\ i\in I\}.
$$
Since each $\beta_i$ is a continuous function, the set $D_I$ is open relative to $Q_I$. The sets $F^\circ_I$ and $D_I$ are mutually orthogonal since the sets $Q_I$ and $H_I$ are by Lemma \ref {L10a}. 

Let ${\bf x}\in D_I\oplus F_I^\circ$ be arbitrary point. Then ${\bf x}={\bf a}+{\bf b}$, where ${\bf a}\in D_I$ and ${\bf b}\in F^\circ _I$. Since $D_I$ and $F^\circ_I$ are relatively open, there is $\epsilon>0$ such that $B({\bf a},\epsilon)\cap Q_I\subset D_I$ and $B({\bf b},\epsilon)\cap H_I\subset F^\circ_I$. By Lemma \ref {L10a}, $Q_I\oplus H_I=\RR^d$. Then for any ${\bf y}\in B({\bf x},\epsilon)$, we have ${\bf y}={\bf q}+{\bf h}$ for some ${\bf q}\in Q_I$ and ${\bf h}\in H_I$. Since $\left|{\bf y}-{\bf x}\right|^2=\left|{\bf q}-{\bf a}\right|^2+\left|{\bf h}-{\bf b}\right|^2<\epsilon^2$, we have $\left|{\bf q}-{\bf a}\right|<\epsilon$ and $\left|{\bf h}-{\bf b}\right|<\epsilon$. Then ${\bf q}\in D_I$ and ${\bf h}\in F^\circ_I$. Then ${\bf y}\in D_I\oplus F^\circ_I$, which implies that $B({\bf x},\epsilon)\subset  D_I\oplus F^\circ_I$. Thus, $ D_I\oplus F^\circ_I$ is open.\hfill $\square$


\begin {thebibliography}{99}
\bibitem {Bab76a}
V.F. Babenko, Asymptotically sharp bounds for the remainder for the best quadrature formulas for several classes of functions, {\it Math. Notes}, {\bf 19} (1976), no. 3, 187--193.
 \bibitem {Bab1978}
 V.F. Babenko, Interpolation of continuous mappings by piecewise
 linear ones, {\it Math. Notes} {\bf 24} (1978), no. 1--2, 526--532.

 \bibitem {Bab1980}
 V.F. Babenko, On optimal piecewise linear interpolation of
 continuous mappings, {\it Investigations in current problems in
 summation and approximation of functions and their applications}
 3--7, Dnepropetrovsk. Gos. Univ., Dnepropetrovsk, 1980 (in Russian).
 
\bibitem {BabBojBor}
V.F. Babenko, B.D. Bojanov, S.V. Borodachov, Optimal recovery of functions on classes defined by a majorant for the modulus of continuity (to appear).
\bibitem {BabBorSko}
V.F. Babenko, S.V. Borodachov, D.S. Skorokhodov, Optimal recovery of isotropic classes of twice differentiable
multivariate functions, {\it Journal of Complexity}, {\bf 26} (2010), 591--607.
\bibitem {BabLyg1975}
 V.F. Babenko, A.A. Ligun, Interpolation by polygonal functions, {\it
 Math. Notes} {\bf 18} (1975), no. 5--6, 1068--1074.
\bibitem {Bah71}
N.S Bakhvalov, On the optimality of linear methods for operator approximation in convex classes of functions,
{\it U.S.S.R. Comput. Math. and Math. Phys.}, {\bf 11} (1971), 244--249.
\bibitem {Boj1975}
 B.D. Bojanov, Best interpolation methods for certain classes of
 differentiable functions, {\it Math. Notes}, {\bf 17} (1975), no.
 3--4, 301--309.


\bibitem{BorSor2010} 
S.V. Borodachov, T.S. Sorokina, An optimal multivariate spline method for recovery of twice differentiable
multivariate functions, {\it BIT Numerical Mathematics} {\bf 51} (2011), no. 3, 497--511.
\bibitem {BorSor2012}
S.V. Borodachov and T.S. Sorokina, Optimal recovery of twice differentiable functions based on symmetric splines, {\it Journal of Approximation Theory} 164 (2012), no. 10, 1443--1459.
\bibitem {CinVanRCA}
E. \c {C}inlar, R.J. Vanderbei, {\it Real and convex analysis}, Undergraduate Texts in Mathematics. Springer, New York, 2013. 
\bibitem {Kif57}
J. Kiefer, Optimum sequential search and approximation methods under minimum regularity assumptions, {J. Soc. Indust. Appl. Math.} {\bf 5} (1957), 105--136. 
\bibitem {Klz1996}
Yu.A. Kilizhekov, Error of approximation by interpolation polynomials of the first degree on $n$-simplices,
{\it Math. Notes} {\bf 60} (1996), no. 3--4, 378--382.
\bibitem {KorTCTP}
 N.P. Korneichuk, {\it Exact constants in approximation theory},
 Cambridge University Press, 1991.
\bibitem {KorSTP}
N. P. Korneichuk, {\it Splines in Approximation Theory}, Nauka, Moscow, 1984 (in Russian).
\bibitem{lawson}
Ch. L. Lawson, Properties of $n$-dimensional triangulations, {\it CAGD} {\bf 3} (1986), 231--246.
\bibitem {Nik50}
S.M. Nikol'skiy, To the question on estimates of approximation by
quadrature formulae, {\it Uspekhi Mat. Nauk}, {\bf 5}, (1950), no.
2(36), 165--177 (in Russian).
\bibitem {NovWozTMP}
E. Novak, H. Wo{\'z}niakowski, {\it Tractability of multivariate problems}, Volumes I, II, and III, European Mathematical Society, (EMS), Z${\rm \ddot {u}}$rich, 2008.
\bibitem {MicRivWin1976}
C.A. Micchelli, T.J. Rivlin, S. Winograd, Optimal recovery of smooth functions, {\it Numer. Math.} {\bf 26} (1976), 191--200.
\bibitem {MotLygDor}
V.P. Motornyi, A.A. Ligun, V.G. Doronin, {\it Optimal recovery of
functions and functionals}, Dnepropetrovsk Univ. Publishers,
Dnepropetrovsk, Ukraine, 1994.
\bibitem {Sar49}
A. Sard, Best approximate integration formulas; best approximation formulas, {\it Amer. J. Math.} {\bf 71} (1949), 80--91.
\bibitem {Smo65}
S.A. Smolyak, On Optimal restoration of functions and functionals of
them. Candidate (Ph.D.) Dissertation, Moscow State University, 1965,
(in Russian).
 \bibitem {TraWozOTOA}
 J.F. Traub, H. Wo{\'z}niakowski, {\it A general theory of optimal
 algorithms,} ACM Monograph Series. Academic Press, Inc. [Harcourt
 Brace Jovanovich, Publishers], New York-London, 1980. xiv+341 pp.
 \bibitem {TraWozWasIBC}
 J.F. Traub, G.W. Wasilkowski, H. Wo{\'z}niakowski, {\it
 Information-based complexity}, Academic Press, Inc., Boston, MA,
 1988. xiv+523 pp.
\bibitem {Zhe03}
A.A. Zhensykbaev, {\it Problems of
recovery of operators}, Institute of computer research,
Moscow-Izhevsk, 2003, 412 pp (in Russian)(ISBN: 5-93972-268-7).
\end {thebibliography}
\end {document}